\documentclass[12pt]{article}
\usepackage{epsfig}
\usepackage{amsthm}
\usepackage{amsfonts}
\usepackage{color}
\font\tenopen=msbm10
\font\sevenopen=msbm7
\font\fiveopen=msbm5
\newfam\openfam
\def\open{\fam\openfam\tenopen}
\textfont\openfam=\tenopen \scriptfont\openfam=\sevenopen
\scriptscriptfont\openfam=\fiveopen \font\title=cmbx10
scaled\magstep1 

\def\N{{\mbox{\open N}}}
\def\Z{{\mbox{\open Z}}}

\def\R{{\mbox{\open R}}}

\def\p{{\mbox{\open P}}}
\def\E{{\mbox{\open E}}}

\font\n=cmcsc10
\def\qed{\hfill$\diamondsuit$}

\theoremstyle{remark} \theoremstyle{lemma}
\theoremstyle{definition} \theoremstyle{corol}
\theoremstyle{proposition} \theoremstyle{condition}
\newtheorem{theorem}{\n{Theorem}}
\newtheorem{remark}{\n{Remark}}
\newtheorem{lemma}{\n{Lemma}}

\newtheorem{corollary}{\n{Corollary}}

\newtheorem{proposition}{\n{Proposition}}

\begin{document}

\centerline{\bf EMPIRICAL PROCESSES OF DEPENDENT RANDOM VARIABLES
\footnote{ 
\newline
{\it $\quad$ Key words and phrases}.  Empirical process, Gaussian
process, Hardy inequality, linear process, Long- and short-range
dependence, Markov chains, martingale, maximal inequality,
nonlinear time series, Sobolev space, tightness, weak
convergence}}

\bigskip

\centerline{\textsc{By Wei Biao Wu}}

\centerline{\today}

\centerline{\it{University of Chicago}}

\bigskip
\pagenumbering{arabic} \noindent {\it Abstract:} Empirical
processes for stationary, causal sequences are considered. We
establish empirical central limit theorems for classes of
indicators of left half lines, absolutely continuous functions and
piecewise differentiable functions. Sample path properties of
empirical distribution functions are also discussed. The results
are applied to linear processes and Markov chains.

\bigskip

\section{Introduction}
The theory of empirical processes plays a fundamental role in
statistics and it has many applications ranging from parameter
estimation to hypothesis testing (van der Vaart and Wellner,
1996). The literature of empirical processes for independent
random variables is huge and there are many deep results; see
Donsker (1952), Dudley (1978), Pollard (1984), Gin\'e and Zinn
(1984), Shorack and Wellner (1986), Ossiander (1987), van der
Vaart and Wellner (1996).

To deal with random variables such as time series that are
dependent, one naturally asks whether results obtained under the
independence assumption remain valid. Such asymptotic theory is
evidently useful for statistical inference of stochastic
processes. Without the independence assumption, it is more
challenging to develop a weak convergence theory for the
associated empirical processes. One way out is to impose strong
mixing conditions to ensure the asymptotic independence; see
Billingsley (1968), Gastwirth and Rubin (1975), Withers (1975),
Mehra and Rao (1975), Doukhan et al (1995), Andrews and Pollard
(1994), Shao and Yu (1996), Rio (1998, 2000) and Pollard (2002)
among others. Other special processes that have been discussed
include linear processes and Gaussian processes; see Dehling and
Taqqu (1989) and Cs\"org\H{o} and Mielniczuk (1996) for long and
short-range dependent subordinated Gaussian processes and Ho and
Hsing (1996) and Wu (2003a) for long-range dependent linear
processes. A collection of recent results is presented in Dehling,
Mikosch and Sorensen (2002). In that collection Dedecker and
Louhichi (2002) made an important generalization of Ossiander's
(1987) result.

Here we investigate the empirical central limit problem for
dependent random variables from another angle that avoids strong
mixing conditions. In particular, we apply a martingale method and
establish a weak convergence theory for stationary, causal
processes. Our results are comparable with the theory for
independent random variables in that the imposed moment conditions
are optimal or almost optimal. We show that, if the process is
short-range dependent in a certain sense, then the limiting
behavior is similar to that of iid random variables in that the
limiting distribution is a Gaussian process and the norming
sequence is $\sqrt n$. For long-range dependent linear processes,
one needs to apply asymptotic expansions to obtain $\sqrt
n$-norming limit theorems (Section \ref{sec:lplrd}).

The paper is structured as follows. In Section \ref{sec:method} we
introduce some mathematical preliminaries necessary for the weak
convergence theory and illustrate the essence of our approach. Two
types of empirical central limit theorems are established.
Empirical processes indexed by indicators of left half lines,
absolutely continuous functions, and piecewise differentiable
functions are discussed in Sections \ref{sec:indicator},
\ref{sec:abc} and \ref{sec:piecewise} respectively. Applications
to linear processes and iterated random functions are made in
Section \ref{sec:cond1}. Section \ref{sec:inequalities} presents
some integral and maximal inequalities that may be of independent
interest. Some proofs are given in Sections \ref{sec:thm12pr} and
\ref{sec:proofThm2}.

\section{Preliminaries}
\label{sec:method} Let $(X_n)_{n \in \Z}$ be a stationary process.
Denote by $F$ and $F_n(x) = n^{-1} \sum_{i=1}^n {\bf 1}_{X_i \le
x}$, $x\in \R$, the marginal and empirical distribution functions.
Let ${\cal G}$ be a class of measurable functions from $\R$ to
$\R$. The centered ${\cal G}$-indexed empirical process is given
by
\begin{eqnarray}
(P_n - P) g = {1\over n} \sum_{i=1}^n g(X_i) - \E [g(X_1)] =
\int_\R g(x) d [F_n(x) - F(x)],\quad g \in {\cal G}.
\end{eqnarray}
The weak convergence theory concerns the limiting behavior of $\{
(P_n - P) g: ~ g \in {\cal G} \}$ under proper scaling. Our
primary goal is to establish abstract Donsker theorems for
dependent random variables. Namely, we aim at finding appropriate
conditions on $(X_n)$ and ${\cal G}$ such that $\{ \sqrt n (P_n -
P) g, ~ g \in {\cal G}\}$ converges in distribution to some tight
Gaussian process $W = \{W(g), ~ g \in {\cal G}\}$.

To describe the weak convergence theory, some mathematical
apparatus is needed. Let the triple $(\Omega, {\cal B}(\Omega),
\p)$ be the probability space on which the process $(X_i)_{i \in
\Z}$ is defined; let $\ell^\infty ({\cal G})$ be the set of
functions $z: {\cal G} \mapsto \R$ for which $\| z \|_{\cal G} : =
\sup_{g \in {\cal G}} |z(g)| < \infty$. A random element $\xi$
with values in $\ell^\infty ({\cal G})$ is said to be tight if,
for every $\delta > 0$, there is a compact set $K_\delta$ of
$(\ell^\infty ({\cal G}), \| \cdot \|_{\cal G})$ for which $\p(\xi
\not \in K_\delta) < \delta$. Assume that for any $x \in \R$,
$\sup_{g \in {\cal G}}|g(x) - P g | < \infty$. Then $\sqrt n (P_n
- P)$ is a map from $\Omega$ to $\ell^\infty ({\cal G})$. We say
that $\sqrt n (P_n - P)$ converges weakly to the tight Gaussian
process $W$ if, for any bounded, continuous function $h:~
(\ell^\infty ({\cal G}), \| \cdot \|_{\cal G}) \mapsto \R$,
$\lim_{n \to \infty} \E^* \{ h[ \sqrt n (P_n - P)]\} = \E [h(W)]$.
Here $\E^*$ denotes the {\it outer expectation}: $\E^* Z = \inf\{
\E U:~ U \ge Z \mbox{ and } U \mbox{ is measurable} \}$. The outer
probability of an arbitrary set $A \subset \Omega$ is given by
$\p^*(A) = \E^* ({\bf 1}_A)$. The outer expectation is introduced
to deal with the measurability issue which emerges when ${\cal G}$
is uncountable. See van der Vaart and Wellner (1996) for more on
weak convergence theory.

It is well-known that $\{ \sqrt n (P_n - P) g, ~ g \in {\cal G}\}
$ converges weakly to some tight Gaussian process $\{W(g), ~ g \in
{\cal G}\}$ if and only if the following two conditions are
satisfied (Theorem 1.5.4, van der Vaart and Wellner, 1996):

(i) {\it Finite-dimensional convergence:} For any finite set of
functions $g_1, \ldots, g_k \in {\cal G}$,
\begin{eqnarray}
\label{eq:finiteclt} [\sqrt n (P_n - P) g_1, \ldots, \sqrt n (P_n
- P) g_k] \Rightarrow \mbox{ a multivariate normal distribution.}
\end{eqnarray}

(ii) {\it Tightness:} For every $\delta, \eta > 0$, ${\cal G}$ can
be partitioned into finitely many sets ${\cal G}_1, \ldots, {\cal
G}_p$ such that
\begin{eqnarray}
\label{eq:tightdef} \limsup_{n\to \infty} \p^* \left[ \max_{1\le
i\le p} \sup_{g, h \in {\cal G}_i} |\sqrt n (P_n - P) (g-h) | \ge
\delta \right] \le \eta.
\end{eqnarray}

The finite-dimensional convergence (i) is trivial if $X_k$ are iid
and $\E[ g^2(X_1)] < \infty$. For dependent random variables, the
central limit theorem (\ref{eq:finiteclt}) itself is an
interesting and important problem and it has received much
attention for more than a half century. Various strong mixing
assumptions are imposed in early work [cf. Rosenblatt (1956),
Ibragimov (1962), Bradley (2002), Doukhan (1994) and Peligrad
(1996) among others]. Gordin (1969) and Gordin and Lifsic (1978)
proposed the martingale approximation method which does not
require strong mixing conditions. This method is quite powerful in
establishing central and non-central limit theorems; see Wu
(2003a, 2003b, 2004a, 2004b), Wu and Woodroofe (2004), Hsing and
Wu (2004) and Wu and Shao (2004), Gordin and Holzmann (2004) for
some recent developments. Earlier contributions can be found in
Hall and Heyde (1980), Kipnis and Varadhan (1986), Voln\'y (1993)
and Ho and Hsing (1997).

Generally speaking, it is more challenging to verify the tightness
condition (ii). Consider the special case of the indicator
function class ${\cal G} = \{ {\bf 1}_{\cdot \le s}:~ s \in \R\}$.
Let $R_n(s) = \sqrt n [F_n(s) - F(s)]$. Then (\ref{eq:tightdef})
requires that there exists a partition $\R = \cup_{i = 1}^p T_i$
such that
\begin{eqnarray*}
\limsup_{n\to \infty} \p \left[ \max_{1\le i\le p} \sup_{s, t \in
T_i} |R_n(s) - R_n(t) | \ge \delta \right] \le \eta.
\end{eqnarray*}
The verification of the preceding relation is not simple partly
due to the fact that $F_n$ is discontinuous. For iid random
variables, there exist sophisticated tools such as chaining
arguments and exponential inequalities. It is not straightforward
to apply those tools to dependent random variables.

In our problem the major interest is to obtain comparable results
without the iid assumption. It is necessary to impose certain
structural assumptions on the underlying process $(X_n)_{n \in
\Z}$ and the class ${\cal G}$. In the early literature, strong
mixing conditions have been imposed on $(X_n)$ [for example,
Doukhan et al (1995), Pollard (2002)]. Here we restrict ourselves
to causal processes. Let $(\varepsilon_k)_{k\in \Z}$ be
independent and identically distributed (iid) random variables;
let $J$ be a measurable function such that
\begin{eqnarray}
\label{eq:Xn} X_n = J(\ldots, \varepsilon_{n-1}, \varepsilon_n)
\end{eqnarray}
is a proper random variable. Such processes are also known as
causal Bernoulli shifts and they have received considerable
attention recently; see Doukhan and Louhichi (1999), Prieur
(2002), Doukhan (2003) among others. The framework (\ref{eq:Xn})
is general enough to allow many interesting and important
examples. Prominent ones are linear processes and Markov chains
arising from iterated random functions; see Section
\ref{sec:cond1}. For the class ${\cal G}$, we consider indicators
of left half lines, absolutely continuous functions and piecewise
differentiable functions.

In this article we will apply the martingale method, and thus shed
new light on this important problem. To illustrate the idea behind
our approach, let ${\cal F}_n = (\ldots, \varepsilon_{n-1},
\varepsilon_{n})$ and denote by $F_\varepsilon( x | {\cal F}_n) =
\p(X_{n+1} \le x | {\cal F}_n)$ the conditional distribution
function of $X_{n+1}$ given the sigma algebra $\sigma({\cal
F}_n)$. Assume throughout the paper that the conditional density
$f_\varepsilon(x | {\cal F}_n) = (\partial / \partial x)
F_\varepsilon(x | {\cal F}_n)$ exists almost surely. Define the
{\it conditional empirical distribution function} ${\tilde F}_n(x)
= n^{-1} \sum_{i=1}^n F_\varepsilon(x | {\cal F}_i)$. Then
\begin{eqnarray}
\label{eq:DecompF} F_n(x) - F(x) = [F_n(x) - {\tilde F}_n(x)] +
[{\tilde F}_n(x) - F(x)].
\end{eqnarray}
The decomposition (\ref{eq:DecompF}) has two important and useful
properties. First, $n[F_n(x) - {\tilde F}_n(x)]$ is a martingale
with stationary, ergodic, and bounded martingale differences.
Second, the function ${\tilde F}_n - F$ is differentiable with
derivative ${\tilde f}_n(x) - f(x)$, where ${\tilde f}_n(x) =
(\partial / {\partial x}) {\tilde F}_n(x) = n^{-1} \sum_{i=1}^n
f_\varepsilon (x | {\cal F}_i)$. These two properties are useful
in establishing tightness. Wu and Mielniczuk (2002) gave a similar
decomposition in the asymptotic theory for kernel density
estimators of linear processes.

The following notation will be used throughout the paper. Let
$w_{\lambda} (d u) = (1+ |u|)^{\lambda}(d u)$ be a weighted
measure. For a random variable $\xi$ write $\xi \in L^q$, $q > 0$,
if $\|\xi \|_q := [\E( |\xi|^q) ]^{ 1/q } < \infty$. Write the
$L^2$ norm $\|\xi\| = \|\xi \|_2$. Define projections ${\cal P}_k
\xi = \E( \xi | {\cal F}_{k}) - \E( \xi | {\cal F}_{k-1})$, $k \in
\Z$. Denote by $C_q$ (resp. $C_\gamma$, $C_\mu$ etc) generic
positive constants which only depend on $q$ (resp. $\gamma$, $\mu$
etc). The values of those constants may vary from line to line.
For a function $p(\cdot)$, let $p {\cal G} = \{p(\cdot) g(\cdot):
~ g \in {\cal G}\}$. For two sequences of real numbers $(a_n)$ and
$(b_n)$, write $a_n \sim b_n$ if $\lim_{n \to \infty} a_n / b_n =
1$.

\section{Empirical distribution functions}
\label{sec:indicator} In this section we consider sample path
properties and weak convergence of empirical distribution
functions. Recall $R_n(s) = \sqrt n [F_n(s) - F(s)]$. The
classical Donsker theorem asserts that, if $X_k$, $k \in \Z$, are
iid random variables, then $\{R_n(s), ~ s\in \R\}$ converges in
distribution to an $F$-Brownian bridge process. The result has
many applications in statistics. To understand the behavior at the
two extremes $s = \pm \infty$, we need to consider the weighted
version $\{ R_n(s) W(s), ~ s\in \R\}$, where $W(s) \to \infty$ as
$s \to \pm \infty$. Clearly, if $W$ is bounded, then by the
continuous mapping theorem, the weak convergence of the weighted
empirical processes follows from that of $R_n$. The
Chibisov-O'Reilly Theorem concerns weighted empirical processes of
iid random variables. A detailed account can be found in Shorack
and Wellner (1986, Section 11.5). The case of dependent random
variables has been far less studied. For strong mixing processes
see Mehra and Rao (1975) and Shao and Yu (1996). Section
\ref{sec:weakedf} generalizes the Chibisov-O'Reilly Theorem to
dependent random variables. Section \ref{sec:modulus} considers
weighted modulus of continuity of $R_n$. Proofs of Theorems
\ref{th:RnW} and \ref{th:modulus} are given in Section
\ref{sec:thm12pr}. Let the weight function $W$ be of the form
$W(t) = (1 + |t|)^\delta$ for some $\delta \ge 0$.

\subsection{A weak convergence result}
\label{sec:weakedf} Let $m$ be a measure on $\R$ and $T_n(\theta)
= \sum_{i=1}^n h(\theta, {\cal F}_i) - n \E [h(\theta, {\cal
F}_1)]$, where $h$ is a measurable function such that $\|
h(\theta, {\cal F}_1) \| < \infty$ for almost all $\theta$ ($m$).
Define
\begin{eqnarray}
\label{eq:sigmahm} \sigma(h, m) = \sum_{j=0}^\infty \sqrt{\int_\R
\| {\cal P}_0 h(\theta, {\cal F}_j) \|^2 m(d \theta) }.
\end{eqnarray}
Let $f'_\varepsilon(\theta | {\cal F}_k) = (\partial / \partial
\theta) f_\varepsilon(\theta | {\cal F}_k)$. In the case that
$h(\theta, {\cal F}_k) = f_\varepsilon(\theta | {\cal F}_k)$ or
$h(\theta, {\cal F}_k) = f'_\varepsilon(\theta | {\cal F}_k)$, we
write $\sigma(f_\varepsilon, m)$ or $\sigma(f'_\varepsilon, m)$
for $\sigma(h, m)$.

\begin{theorem}
\label{th:RnW} Let $\gamma > 0$ and $q > 2$. Assume $\E(
|X_1|^\gamma)< \infty$ and
\begin{eqnarray}
\label{eq:IntfEq} \int_{\R} \E [f_\varepsilon^{q/2} (u | {\cal
F}_0)] w_{\gamma -1 + q/2} (d u) < \infty.
\end{eqnarray}
In addition assume
\begin{eqnarray}
\label{eq:fwep6}
 \sigma(f_\varepsilon, w_{1 + 2\gamma / q})
 +\sigma(f'_\varepsilon, w_{-1 - 2\gamma / q}) < \infty.
\end{eqnarray}
Then (i)
\begin{eqnarray}
\label{eq:Enr}
 \E \left[\sup_{s\in \R} |R_n(s)|^2 (1+|s|)^{2\gamma / q} \right]
 = {\cal O}(1)
\end{eqnarray}
and (ii) the process $\{ R_n(s) (1+|s|)^{\gamma /q}, ~ s\in \R\}$
converges weakly to a tight Gaussian process.
\end{theorem}

An important issue in applying Theorem \ref{th:RnW} is the
verification of (\ref{eq:fwep6}), which is basically a short-range
dependence condition. For many important models such as linear
processes and Markov chains, (\ref{eq:fwep6}) is easily
verifiable; see Section \ref{sec:cond1}. In particular, if $X_n$
is a linear process, then (\ref{eq:fwep6}) reduces to the
conventional definition of the short-range dependence of linear
processes.

\begin{remark}{\rm
Let $\alpha_n = \sup\{ |\p(A \cap B) - \p(A) \p(B)|: A \in {\cal
A}_0, B \in {\cal B}_n\}$ be the strong mixing coefficients, where
${\cal A}_k= \sigma(X_i, i\le k)$ and ${\cal B}_n = \sigma(X_j,
j\ge n)$. If $X_n$ is strong mixing, namely $\alpha_n \to 0$, Rio
(2000) showed that $\E [ \sup_{s \in \R} |R_n(s)|^2] \le (1 + 4
\sum_{k=0}^{n-1} \alpha_k) (2 + \log n)^2$; see Proposition 7.1
therein. Clearly the bound in (\ref{eq:Enr}) is sharper. \qed }
\end{remark}

\begin{corollary}
\label{cor:EmpWeighted} Let $\gamma \ge 0$. Assume that $\E
(|X_1|^{\gamma + \delta}) < \infty$ and $\sup_u f_\varepsilon(u |
{\cal F}_0) \le \tau$ hold for some $\delta > 0$ and $\tau <
\infty$. Further assume that $\sigma(f_\varepsilon, w_{1 + \gamma}
) < \infty$ and $\sigma(f'_\varepsilon, w_{-1 - \gamma}) <
\infty$. Then $\{ R_n(s) ( 1 + |s|)^{ \gamma/2 }, ~ s \in \R \}$
converges weakly to a tight Gaussian process.
\end{corollary}
\noindent{\it Proof}. Let $q = 2(\gamma + \delta + 1)/(\gamma +
1)$ and $\gamma' = \gamma q /2$. Then $\gamma' -1 + q/2= \gamma +
\delta$. We shall apply Theorem \ref{th:RnW} with $\gamma'$ and
$q$. Since $\E [f_\varepsilon (u | {\cal F}_0)] = f(u)$,
(\ref{eq:IntfEq}) follows from
\begin{eqnarray*}
\int_{\R} \E [f_\varepsilon^{q/2} (u | {\cal F}_0)] w_{\gamma' -1
+ q/2} (d u)
 &\le& \tau^{q/2 -1} \int_{\R} \E [f_\varepsilon (u |
{\cal F}_0)] w_{\gamma' -1 + q/2} (d u) \\
 &=& \tau^{q/2 -1} \E[ (1+ |X_1|)^{\gamma + \delta}]
 < \infty.
\end{eqnarray*}
Note that $2 \gamma' / q = \gamma$. Then (\ref{eq:fwep6}) holds
and the corollary follows from Theorem \ref{th:RnW}. \qed

\bigskip

It is interesting to compare Corollary \ref{cor:EmpWeighted} with
the Chibisov-O'Reilly Theorem, which concerns weighted empirical
processes of iid random variables. The moment condition $\E
(|X_1|^{\gamma + \delta}) < \infty$ of Corollary
\ref{cor:EmpWeighted} is almost necessary in the sense that it
cannot be replaced by the weaker one
\begin{eqnarray}
\label{eq:weakempmom} \E\{ |X_1|^{\gamma} \log^{-1}( 2 + |X_1|)
[\log \log (10 + |X_1|)]^{-\lambda} \} < \infty
\end{eqnarray}
for some $\lambda > 0$. To see this let $X_k$ be iid symmetric
random variables with continuous, strictly increasing distribution
function $F$; let $F^\#$ be the quantile function and $m(u) =
[1+|F^\#(u)| ]^{-\gamma/2}$. Then we have the distributional
equality
\begin{eqnarray*}
 \{ R_n(s)(1+|s|)^{\gamma /2}, ~ s\in \R\}
 =_{\cal D} \{ R_n(F^\#(u)) / m(u), ~ u\in (0,1)\},
\end{eqnarray*}
Assume that $F(s) (1+|s|)^\gamma$ is increasing on $(-\infty, G)$
for some $G < 0$. Then $m(u) / \sqrt u$ is decreasing on $(0,
F(G))$. By the Chibisov-O'Reilly Theorem, $\{ R_n(F^\#(u)) / m(u),
u\in (0,1)\}$ is tight if and only if $\lim_{t \downarrow 0} {{
m(t)} / {\sqrt{t \log \log (t^{-1})}}} = \infty$, namely
\begin{eqnarray}
\label{eq:Chibisov}
 \lim_{u \to -\infty} F(u) (1+|u|)^\gamma \log\log |u| = 0.
\end{eqnarray}
The above condition controls the heaviness of the tail of $X_1$.
Let $F(u) = |u|^{-\gamma} (\log \log |u| )^{-1}$ for $u \le -10$.
Then it is easily seen that (\ref{eq:weakempmom}) holds, while
(\ref{eq:Chibisov}) is violated. It is unclear whether stronger
versions of (\ref{eq:weakempmom}) such as $\E( |X_1|^\gamma ) <
\infty$ or $\E[ |X_1|^\gamma \log^{-1}( 2 + |X_1|)] < \infty$ are
sufficient.

\subsection{Modulus of continuity}
\label{sec:modulus} Theorem \ref{th:modulus} below concerns the
weighted modulus of continuity of $R_n(\cdot)$. Sample path
properties of empirical distribution functions of iid random
variables have been extensively explored; see for example
Cs\"org\H o {\it et al.} (1986), Shorack and Wellner (1986) and
Einmahl and Mason (1988) among others. It is far less studied for
the dependent case.

\begin{theorem}
\label{th:modulus} Assume $\E(|X_1|^\gamma)< \infty$ and $\sup_u
f_\varepsilon(u | {\cal F}_0) \le \tau$ for some $\gamma \ge 0$
and $\tau < \infty$. Let $\delta_n < 1/2$ be a sequence of
positive numbers such that $(\log n)^{2 q/(q-2)} = {\cal O}(n
\delta_n)$, where $2 < q < 4$. Further assume that, for some $\mu
\le 1$,
\begin{eqnarray}
\label{eq:fwep}
 \sigma(f_\varepsilon, w_{2\gamma/ q - \mu})
 + \sigma(f'_\varepsilon, w_{2\gamma/ q + \mu}) < \infty.
\end{eqnarray}
Then there exists a constant $0 < C < \infty$, independent of $n$
and $\delta_n$, such that for all $n \ge 1$,
\begin{eqnarray}
\label{eq:wmodulus}
 \E \left[\sup_{t\in \R} ( 1 + |t|)^{2\gamma / q}
 \sup_{|s| \le \delta_n} |R_n(t+s) - R_n(t)|^2 \right]
 \le C \delta_n^{1 - 2/q}.
\end{eqnarray}
\end{theorem}

\section{Absolutely continuous functions}
\label{sec:abc} Let ${\cal AC}$ be the collection of all
absolutely continuous functions $g: \R\mapsto \R$. In this section
we shall investigate the behavior of the empirical process ${\sqrt
n} (P_n - P) g$ indexed by $g \in {\cal G}_{\gamma, \mu}$ or $g
\in {\cal H}_{\eta, \delta}$, where ${\cal G}_{\gamma, \mu}$,
$\gamma \ge 0$, $\mu \le 1$, is the weighted Sobolev class
\begin{eqnarray}
\label{eq:Ggamma} {\cal G}_{\gamma, \mu}
  = \left \{ g \in {\cal AC}:
  ~ \int_\R g^2(u) w_{-\gamma - \mu}(d u)
  + \int_\R [g'(u)]^2 w_{-\gamma + \mu}(d u)
  \le 1 \right\}
\end{eqnarray}
and the class ${\cal H}_{\eta, \delta}$, $\eta, \delta \ge 0$,
$\eta - \delta \le 1$, is given by
\begin{eqnarray}
\label{eq:Heta} {\cal H}_{\eta, \delta} = \{g \in {\cal AC}:~
|g(u)| \le (1+|u|)^\eta,~ |g'(u)| \le (1+|u|)^\delta\mbox{ for all
} u \in \R \}.
\end{eqnarray}
Functions in ${\cal G}_{\gamma, \mu}$ can be unbounded if $\gamma
+ \mu > 1$. See Remark \ref{rem:GH2} for some properties of ${\cal
G}_{\gamma, \mu}$ and ${\cal H}_{\eta, \delta}$. In the study of
empirical central limit theorems for independent random variables,
bracketing conditions are often imposed. The class of
differentiable functions is an important case that bracketing
conditions can be verified; see Chapter 2.7 in van der Vaart and
Wellner (1996). Functions considered in robust inference are often
absolutely continuous.

\begin{theorem}
\label{th:Pn-Pg} (i) Let $\gamma \ge 0$ and $\mu \le 1$. Assume
$\E ( |X_1|^\gamma) < \infty$ and
\begin{eqnarray}\label{eq:srdG}
\sigma(f_\varepsilon, w_{\gamma + \mu}) < \infty.
\end{eqnarray}
Then $\{\sqrt n (P_n - P) g: ~ g \in {\cal G}_{\gamma, \mu}\}$
converges weakly to a tight Gaussian process $\{W(g): ~ g \in
{\cal G}_{\gamma, \mu}\}$ with mean $0$ and finite covariance
function
\begin{eqnarray*}
{\rm cov}[W(g_1), W(g_2)] = \sum_{k \in \Z} {\rm cov}[g_1(X_0),
g_2(X_k)], \quad g_1, g_2 \in {\cal G}_{\gamma, \mu}.
\end{eqnarray*}
(ii) Let $\eta, \delta \ge 0$. Assume that
\begin{eqnarray}
\label{eq:Ginezinn1} \sum_{j=1}^\infty j^{\alpha \eta} \sqrt{ \p[
j^\alpha \le |X_1| < (j+1)^\alpha] } < \infty
\end{eqnarray}
if $\eta - \delta < 1$, where $\alpha = (1 + \delta - \eta)^{-1}$,
and that
\begin{eqnarray}
\label{eq:Ginezinn2} \sum_{j=1}^\infty 2^{j \eta} \sqrt{ \p( 2^j
\le |X_1| < 2^{j+1} ) } < \infty
\end{eqnarray}
if $\eta - \delta = 1$. In addition, let $m(d x) = (1 + |x|)^{1 +
2\eta} \log^2 (2+ |x|) d x$ and assume
\begin{eqnarray}\label{eq:srdH}
\sigma(f_\varepsilon, m) < \infty.
\end{eqnarray}
Then the conclusion in (i) still holds with ${\cal G}_{\gamma,
\mu}$ replaced by ${\cal H}_{\eta, \delta}$.
\end{theorem}

The proof of Theorem \ref{th:Pn-Pg} is given in Section
\ref{sec:thm3pf}. Theorem \ref{th:Pn-Pg} generalizes the empirical
central limit theorems in Gin\'e and Zinn (1986) in two ways: by
allowing dependence and by considering wider classes. Gin\'e and
Zinn (1986) considered iid random variables and Lipschitz
continuous functions. In particular, they show that (i) if $\eta =
\delta =0$ (the class of bounded Lipschitz continuous functions),
then (\ref{eq:Ginezinn1}) is necessary and sufficient for the
tightness of $\{\sqrt n (P_n - P) g:~ g \in {\cal H}_{0, 0}\}$ and
(ii) if $\delta = 0$, then (\ref{eq:Ginezinn2}) is necessary and
sufficient for the tightness of $\{\sqrt n (P_n - P) g: ~ g \in
{\cal H}_{0, \cdot}\}$, where ${\cal H}_{0, \cdot} = \{g: ~
|g'(u)| \le 1 \}$ is the class of Lipschitz continuous functions.
Further consideration is given in van der Vaart (1996). So in the
cases of $\eta = \delta =0$ and $\delta = 0$, $\eta = 1$, the
conditions (\ref{eq:Ginezinn1}) and (\ref{eq:Ginezinn2}) are
optimal.

The dependence structure of $(X_k)_{k \in \Z}$ will inevitably
find its way into the weak convergence theory. Here the assumption
on the dependence is encapsulated as $\sigma(h, m) < \infty$,
which is a consequence of the requirement of $\sqrt
n$-normalization of the partial sums in view of $\int_\R
\|T_n(\theta)\|^2 m(d \theta) \le n \sigma^2(h, m)$ by Lemma
\ref{lem:Intsum}. Conditions (\ref{eq:srdG}) and (\ref{eq:srdH})
actually imply more. They are natural conditions for the
asymptotic normality of $\sqrt n (P_n - P) g$ for a fixed $g \in
{\cal G}_{\gamma, \mu}$ or $g \in {\cal H}_{\eta, \delta}$; see
the proof of Theorem \ref{th:Pn-Pg} in Section \ref{sec:thm3pf}.
If $X_n$ is a long-range dependent linear process, then
(\ref{eq:srdG}) and (\ref{eq:srdH}) are violated and the norming
sequence of $(P_n - P) g$ is different from $\sqrt n$. In this
case asymptotic expansions are needed to ensure a Gaussian limit
process with a $\sqrt n$-normalization (cf. Section
\ref{sec:lplrd}).

Doukhan et al (1995), Rio (1998) and Pollard (2002) considered
stationary, absolutely regular processes. Bracketing conditions
are given based on a metric which involves mixing coefficients.
Rio (2000, Theorem 8.1) established an empirical central limit
theorem for strong mixing processes indexed by Lipschitz
continuous functions. In the case of causal processes, it is not
easy to verify that they are absolutely regular or strong mixing.
For linear processes to be strong mixing, quite restrictive
conditions are needed on the decay rate of the coefficients
(Doukhan, 1994). In comparison, Theorem \ref{th:lp} and Corollary
\ref{cor:EmpWtdlp} (cf. Section \ref{sec:lpsrd}) impose a sharp
condition on the coefficients of linear processes.

Some new weak dependence conditions are used in Doukhan and
Louhichi (1999), Prieur (2002) and Dedecker and Prieur (2003b).
Here we compare our results with theirs by applying them to the
Gaussian process $X_n = \sum_{i=0}^\infty a_i \varepsilon_{n-i}$,
where $a_i = (i+1)^{-\beta}$, $\beta > 1$, and $\varepsilon_n$ are
iid standard normal. As mentioned in the preceding paragraph,
Theorem \ref{th:lp} and Corollary \ref{cor:EmpWtdlp} impose the
sharp condition $\beta > 1$. Proposition 2 in Doukhan and Louhichi
(1999) asserts that $R_n$ is tight if
\begin{eqnarray}\label{eq:dkl99}
\sup_{f \in {\cal I}} |{\rm cov}[f(X_{t_1}) f(X_{t_2}), \,
 f(X_{t_3}) f(X_{t_4})]| = O(n^{-5/2 - \delta})
\end{eqnarray}
for some $\delta > 0$, where ${\cal I} = \{ x \to {\bf 1}_{s < x <
t}:~ s, t \in \R\}$, $0 \le t_1 \le t_2 \le t_3 \le t_4$ and $n =
t_3 - t_2$.  Let $t_1 = t_2 = 0$, $t_3 = t_4 = n$ and $f(x) = {\bf
1}_{x < 0}$. Elementary manipulations show that the covariance
$r(n) = \E (X_0 X_n) \sim C_\beta n^{1- 2 \beta}$ for some
constant $0 < C_\beta < \infty$ and $| \p(X_0 < 0, \, X_n < 0) -
\p(X_0 < 0) \p(X_n < 0)| = O[r(n)]$. So (\ref{eq:dkl99}) requires
$1- 2 \beta \le -5/2 - \delta$, namely $\beta > 7/4$. The
empirical central limit theorem in Dedecker and Prieur (2003b) is
not applicable to the process $X_n$. Actually, Corollary 4 in the
latter paper assumes that $\sum_{k=1}^\infty \phi(k) < \infty$,
where $\phi(k) = \sup_{t\in \R} \| \p (X_k \le t | {\cal F}_0) -
\p (X_k \le t ) \|_\infty$. It is easily seen that $\phi(k) \equiv
1$ for the process $X_k$. By Theorem 1 in Prieur (2002), $R_n$
converges weakly if $\E| X_n - \sum_{i=0}^n a_i \varepsilon_{n-i}|
= O(n^{-2 - 2\sqrt 2 - \delta})$ for some $\delta > 0$. The latter
condition implies $\beta > 2.5 + 2 \sqrt 2$. On the other hand,
the existence of the conditional density $f_\varepsilon(x | {\cal
F}_n)$ is not assumed in those papers, while it is needed in our
results.

\begin{remark}
\label{rem:GH2} {\rm We say that two classes ${\cal G}_1$ and
${\cal G}_2$ are equivalent, denoted by ${\cal G}_1 \sim {\cal
G}_2$, if there is a constant $\lambda$ such that $\lambda^{-1}
{\cal G}_1 \subset {\cal G}_2 \subset \lambda {\cal G}_1$. In this
case the process $\{\sqrt n (P_n - P) g:~ g \in {\cal G}_1\}$
converges weakly if and only if $\{\sqrt n (P_n - P) g:~ g \in
{\cal G}_2\}$ does.

If $\eta \ge 1+ \delta$, then ${\cal H}_{\eta, \delta} \sim {\cal
H}_{\delta + 1, \delta}$. Clearly ${\cal H}_{\delta + 1, \delta}
\subset {\cal H}_{\eta, \delta}$. Let $u \ge 0$. Then
\begin{eqnarray*}
|g(u)| \le |g(0)| + \int_0^u |g'(t)| d t \le 1 + \int_0^u (1 +
t)^\delta d t \le (1 + u)^{1 + \delta}
\end{eqnarray*}
and hence ${\cal H}_{\eta, \delta} \subset {\cal H}_{\delta + 1,
\delta}$. The two classes ${\cal G}_{\gamma, \mu}$ and ${\cal
H}_{\eta, \delta}$ are closely related. If $\gamma > 1 + \eta +
\delta$, then ${\cal H}_{\eta, \delta} \subset C_{\gamma, \eta,
\delta} {\cal G}_{\gamma, \eta - \delta}$ for some $C_{\gamma,
\eta, \delta} < \infty$. There is no converse inclusion.

A particularly interesting case is when $\gamma > 0$ and $\mu =
1$. Let
\begin{eqnarray}
\label{eq:Gstar} {\cal G}_{\gamma, \mu, 0}^* = \{ g \in {\cal AC}:
~ g(0) = 0 \mbox{ and } \int_\R [g'(u)]^2 w_{-\gamma + \mu}(d u)
\le 1 \}
\end{eqnarray}
and ${\cal G}_{\gamma, \mu, 0} = \{g \in {\cal G}_{\gamma, \mu}: ~
g(0) = 0\}$. By Hardy's inequality (\ref{eq:IntHmu1}), ${\cal
G}_{\gamma, 1, 0} \subset (1 + 4\gamma^{-2})^{1/2} {\cal
G}_{\gamma, 1, 0}^*$. The other relation ${\cal G}_{\gamma, 1,
0}^* \subset {\cal G}_{\gamma, 1, 0}$ is obvious. Therefore ${\cal
G}_{\gamma, 1, 0}^* \sim {\cal G}_{\gamma, 1, 0}$. By Lemma
\ref{lem:supHweighted}, for $g \in {\cal G}_{\gamma, \mu}$,
$\sup_x g^2(x)(1 + |x|)^{-\gamma} \le C_{\gamma, \mu}$ for some
$C_{\gamma, \mu} < \infty$. So $g^2(0) \le C_{\gamma, \mu}$ and
there exists a constant $C_{\gamma, \mu}'$ such that $g - g(0) \in
C_{\gamma, \mu}' {\cal G}_{\gamma, \mu}$ whenever $g \in {\cal
G}_{\gamma, \mu}$. Therefore ${\cal G}_{\gamma, 1, 0} \sim \{g -
g(0): ~ g \in {\cal G}_{\gamma, 1} \}$. Notice that $(P_n - P) (g
- g(0)) = (P_n - P) g$. So the weak convergence problem of
$\{\sqrt n (P_n - P) g:~ g \in {\cal G}_{\gamma, 1}\}$ is
equivalent to the seemingly simpler one $\{\sqrt n (P_n - P) g: ~
g \in {\cal G}_{\gamma, 1, 0}^*\}$. } \qed
\end{remark}

\subsection{Proof of Theorem \ref{th:Pn-Pg}.} \label{sec:thm3pf}
We first consider the case ${\cal G} = {\cal G}_{\gamma, \mu}$.
Following the ideas behind the decomposition (\ref{eq:DecompF}),
write $\sqrt n (P_n - P) g = M_n(g) + N_n(g)$, where
\begin{eqnarray*}
M_n(g) = n^{1/2} \int_\R g(x) d [ F_n(x) - {\tilde F}_n (x) ]
\mbox{ and } N_n(g) = n^{1/2} \int_\R g(x) [{\tilde f}_n (x) -
f(x)] d x.
\end{eqnarray*}
The tightness of the process $\{\sqrt n (P_n - P) g:~ g \in {\cal
G}_{\gamma, \mu}\}$ follows from that of $\{M_n(g):~ g \in {\cal
G}_{\gamma, \mu}\}$ and $\{N_n(g):~ g \in {\cal G}_{\gamma,
\mu}\}$, which are asserted by Propositions \ref{prop:Mngtight}
and \ref{prop:Nngtight} (cf. Sections \ref{sec:Mn} and
\ref{sec:Nn}) respectively. It remains to establish the
finite-dimensional convergence. By Lemma \ref{lem:supHweighted},
there exists $C_{\gamma, \mu} < \infty$ such that
\begin{eqnarray*}
 \sup_{x \in \R}  { {g^2(x)} \over { (1+|x|)^{\gamma} } }
 \le C_{\gamma, \mu} \int_\R g^2(u) w_{- \gamma - \mu} (d u)
 + C_{\gamma, \mu} \int_\R [g'(u)]^2 w_{- \gamma + \mu} (d u)
 \le C_{\gamma, \mu}.
\end{eqnarray*}
Hence $\| g(X_1)\|^2 \le C_{\gamma, \mu} \E[(|X_1|+1)^\gamma] <
\infty$. Notice that for $n \ge 1$,
\begin{eqnarray*}
{\cal P}_0 g(X_n) = {\cal P}_0 \E[g(X_n) | {\cal F}_{n-1}]
 = {\cal P}_0 \int_\R g(x) f_\varepsilon(x| {\cal F}_{n-1}) d x
 = \int_\R g(x) {\cal P}_0 f_\varepsilon(x| {\cal F}_{n-1}) d x.
\end{eqnarray*}
By the Cauchy-Schwarz inequality,
\begin{eqnarray*}
\|{\cal P}_0 g(X_n)\| \le \sqrt{\int_\R g^2(x) w_{-\gamma - \mu}(d
x) \int_\R \| {\cal P}_0 f_\varepsilon (x| {\cal F}_{n-1})\|^2
w_{\gamma + \mu}(d x)},
\end{eqnarray*}
which in view of (\ref{eq:srdG}) is summable. By Lemma
\ref{lem:hannan}, the finite-dimensional convergence holds.

The other case that ${\cal G} = {\cal H}_{\eta, \delta}$ similarly
follows. For $g \in {\cal H}_{\eta, \delta}$, $|g(u)| \le (1 +
|u|)^\eta$. By (\ref{eq:Ginezinn1}) or (\ref{eq:Ginezinn2}), we
have $\E(|X_1|^{2 \eta}) < \infty$, which by Lemma
\ref{lem:hannan} implies the finite-dimensional convergence. The
tightness follows from (ii) of Proposition \ref{prop:Mngtight} and
Proposition \ref{prop:Nngtight}. \qed

\subsection{Tightness of $M_n$}
\label{sec:Mn} In this section we establish the weak convergence
of $M_n$. The tightness of $M_n$ involves moment conditions on
$X_1$ and sizes of ${\cal G}_{\gamma, \mu}$ or ${\cal H}_{\eta,
\delta}$, which are characterized by the parameters $(\gamma,
\mu)$ and $(\eta, \delta)$. There is a tradeoff between the moment
conditions and the sizes of the classes: larger classes require
stronger moment conditions.

Interestingly, the tightness of $M_n$ does not involve the
dependence structure of $(X_k)_{k \in \Z}$. In Section
\ref{sec:lplrd}, we apply the result to strongly dependent
processes.

\begin{proposition}
\label{prop:Mngtight}  (i) Assume that $\E(|X_1|^\gamma) <
\infty$. Then the process $\{M_n(g):~ g \in {\cal G}_{\gamma,
\mu}\}$ is tight. Consequently, it converges weakly to a tight
Gaussian process $\{W_M(g): ~ g \in {\cal G}_{\gamma, \mu}\}$ with
mean $0$ and covariance function
\begin{eqnarray}
\label{eq:covWM} {\rm cov}[W_M(g_1), W_M(g_2)] = \E\{ [ {\cal P}_0
g_1(X_1)][ {\cal P}_0 g_2(X_1)] \}, \quad g_1, g_2 \in {\cal
G}_{\gamma, \mu}.
\end{eqnarray}
(ii) Let $\eta, \delta \ge 0$. Assume (\ref{eq:Ginezinn1}) if
$\eta - \delta < 1$ and (\ref{eq:Ginezinn2}) if $\eta - \delta =
1$. Then $\{M_n(g): ~ g \in {\cal H}_{\eta, \delta} \}$ converges
weakly to a centered, tight Gaussian process $\{W_M(g): ~ g \in
{\cal H}_{\eta, \delta}\}$ with the covariance function
(\ref{eq:covWM}).
\end{proposition}

The proof of Proposition \ref{prop:Mngtight} is given in Section
\ref{sec:proofThm2}. Dedecker and Louhichi (2002) generalized
Ossiander's (1987) method to dependent sequences. They apply their
results to martingales (cf. Theorem 3.3 and Section 4.2 therein).
A related result is given in Nishiyama (2000). The conditions
imposed in the latter paper seem hard to work with; see van der
Vaart and Wellner (1996, p. 212). Here we compare our Proposition
\ref{prop:Mngtight} with that of Dedecker and Louhichi (2002).

For a class ${\cal G}$ and a norm $d$ let $N(\delta, {\cal G}, d)$
be the minimum number of $\delta$-brackets needed to cover ${\cal
G}$. Here for two functions $l$ and $u$, the bracket $[l,u]$ is
defined by $[l,u] = \{g \in {\cal G}: l(x) \le g(x) \le u(x)
\mbox{ for all } x\}$. If $d(u-l) \le \delta$, then we say that
$[l,u]$ is a $\delta$-bracket. For a function $g$ define its
essential supremum $d_2$ norm by
\begin{eqnarray}
\label{eq:Dedd2} d^2_2(g) = {\rm ess ~ sup} \E [ g^2(X_1) | {\cal
F}_0].
\end{eqnarray}
Dedecker and Louhichi (2002) show that the $\ell^\infty$-valued
map $\sqrt n(P_n - P)$ is asymptotically $d_2$-equicontinuous if
\begin{eqnarray}
\label{eq:Dedd3}
 \int_0^1 \sqrt{ \log N(u, {\cal G}, d_2)} d u < \infty.
\end{eqnarray}
However, the norm $d_2$ in (\ref{eq:Dedd2}) is so strong that
(\ref{eq:Dedd3}) is violated for many important applications. For
example, let $X_n = \sum_{i=0}^\infty a_i \varepsilon_{n - i}$,
where $\varepsilon_k$ are iid random variables. Assume that the
support of $\varepsilon_1$ is the whole real line, $a_0 = 1$ and
that there are infinitely many non-zero coefficients. Let ${\cal
G}_H = \{g(\cdot - \theta), ~ \theta \in \R\}$, where $g(x) =
\max[-1, \min(x, 1)]$ is the derivative of a Huber function.
Elementary calculations show that, since $Y_0 = X_1 -
\varepsilon_1$ has the support of the whole real line, the bracket
number $N(\delta, {\cal G}_H, d_2) = \infty$ for all $0 < \delta <
1/2$. Nevertheless, their generalization is quite important since
the condition (\ref{eq:Dedd3}) is tractable in many cases. We
apply their result in the proof of Proposition \ref{prop:Mngtight}
and overcome the limitation by using a truncation argument.

\begin{remark}{\rm
Let ${\cal G}_{D L}$ be a class of functions such that it
satisfies Dedecker and Louhichi's bracketing condition
(\ref{eq:Dedd3}) and $({\cal G}_{D L}, d_2)$ is totally bounded.
Then it is easily seen that the pairwise sum ${\cal G}_{D L}
+{\cal G}_{\gamma, \mu} = \{ g + h:~ g \in {\cal G}_{D L}, ~h \in
{\cal G}_{\gamma, \mu} \}$ is also a Donsker class, namely
Proposition \ref{prop:Mngtight} holds for this pairwise sum. \qed}
\end{remark}

\subsection{Tightness of $N_n$}
\label{sec:Nn} The tightness of $N_n$ requires the short-range
dependence condition (\ref{eq:srdG}) or (\ref{eq:srdH}). The fact
that ${\tilde F}_n$ is differentiable is quite useful.

\begin{proposition}
\label{prop:Nngtight} Assume (\ref{eq:srdG}) (resp.
(\ref{eq:srdH})). Then the process $\{N_n(g):~ g \in {\cal
G}_{\gamma, \mu}\}$ (resp. $\{N_n(g):~ g \in {\cal H}_{\eta,
\delta}\} $) is tight.
\end{proposition}
\noindent{\it Proof.} Recall ${\bf 1}_{|\cdot| \le r} {\cal
G}_{\gamma, \mu} = \{ g {\bf 1}_{|\cdot| \le r}: ~ g \in {\cal
G}_{\gamma, \mu}\}$. For $g \in {\cal G}_{\gamma, \mu}$,
\begin{eqnarray*}
\left[ \int_{|x| \ge r} | g(x) [ {\tilde f}_n(x) - f(x)] | d x
\right]^2
 &\le& \int_{|x| \ge r} g^2(x) w_{-\gamma - \mu}(d x)
 \int_{|x| \ge r} |{\tilde f}_n(x) - f(x)|^2 w_{\gamma + \mu}(d x)\\
 &\le& \int_{|x|\ge r}|{\tilde f}_n(x)- f(x)|^2 w_{\gamma+\mu}(d x).
\end{eqnarray*}
Applying Lemma \ref{lem:Intsum} with $A= (-\infty, -r) \cup (r,
\infty)$ and $T_n(\theta) = \sum_{i=1}^n f_\varepsilon(\theta
|{\cal F}_j)$, we have
\begin{eqnarray*}
&&\E^* \left\{ \sup_{g \in {\cal G}_{\gamma, \mu}} \left[ \sqrt n
\int_{|x| > r} | g(x) [ {\tilde f}_n(x) - f(x)] | d x \right]^2
\right\} \cr
 &&\le n \E \left\{ \int_{|x| > r} |{\tilde f}_n(x) - f(x)|^2
 w_{\gamma + \mu}(d x) \right\} \\
 && \le \left[ \sum_{j=0}^\infty \sqrt{ \int_{|x| > r}
 \| {\cal P}_0 f_\varepsilon(\theta| {\cal F}_j) \|^2
 w_{\gamma + \mu} (d \theta)} \right]^2
\end{eqnarray*}
which converges to $0$ as $r \to \infty$. As in the proof of
Proposition \ref{prop:Mngtight} (cf. Section \ref{sec:proofThm2}),
it remains to show that for fixed $r$, $\{N_n(g):~ g \in {\bf
1}_{|\cdot | \le r} {\cal G}_{\gamma, \mu}\}$ is tight. To this
end, we apply Lemma~\ref{lem:Intsum} with $A = \R$ and
$T_n(\theta) = \sum_{i=1}^n f_\varepsilon (\theta|{\cal F}_j)$. By
(\ref{eq:srdG}), there exists $\kappa < \infty$ such that
\begin{eqnarray*}
n \int_\R \| {\tilde f}_n(x) - f(x)\|^2 w_{\gamma + \mu}(d x)
 \le \kappa
\end{eqnarray*}
holds for all $n \in \N$. So the asymptotic equi-continuity
follows from
\begin{eqnarray*}
&&\lim_{\delta \to 0} \limsup_{n \to \infty} \p^* \left\{ \sup_{\|
 g-h \|_\infty \le \delta, ~ g, h \in {\bf 1}_{|\cdot| \le r} {\cal
 G}_{\gamma, \mu} } |N_n(g-h)| > \epsilon \right\} \\
 & & \le \lim_{\delta \to 0} \limsup_{n \to \infty} \p \left\{
 \sqrt n \int_{-r}^r | {\tilde f}_n(x) - f(x)| d x >
 {\epsilon \over \delta} \right\} \\
 & & \le \lim_{\delta \to 0} \limsup_{n \to \infty}
 {{\delta^2} \over {\epsilon^2}} n \E
 \left[\int_{-r}^r | {\tilde f}_n(x) - f(x)| w_{\gamma + \mu}(d x)
 \right]^2 =0.
\end{eqnarray*}
It is easily seen that $({\bf 1}_{|\cdot | \le r} {\cal
G}_{\gamma, \mu}, \|\cdot\|_\infty)$ is totally bounded [cf.
(\ref{eq:Sobolev}) and (\ref{eq:GS})]. Therefore, the process
$\{N_n(g):~ g \in {\cal G}_{\gamma, \mu}\}$ is tight.

The other case that $g \in {\cal H}_{\eta, \delta}$ can be
similarly proved by noticing that
\begin{eqnarray*}
& &\sup_{g \in {\cal H}_{\eta, \delta}} \int_A | g(x) [ {\tilde
f}_n(x) - f(x)] | d x
 \le  \int_A | {\tilde f}_n(x) - f(x)|
w _\eta(d x) \\
 && \le \left[ \int_A | {\tilde f}_n(x) - f(x)|^2 m(d
x) \right]^{1/2} \times \left[ \int_A  (1 + |x|)^{-1} \log^{-2}
(2+ |x|) d x \right]^{1/2}
\end{eqnarray*}
in view of the Cauchy-Schwarz inequality. \qed

\subsection{Finite-dimensional convergence}
\begin{lemma}
\label{lem:hannan} Let $g_1, \ldots, g_k$ be measurable functions
such that $g(X_i) \in L^2$ and $\E [g_i(X_1)] = 0$, $1\le i\le k$.
Let ${\bf g} = (g_1, \ldots, g_k)$. Assume
\begin{eqnarray}
\label{eq:bernoulli}
  \sum_{m=0}^\infty \|{\cal P}_0 g_i(X_m) \| < \infty,
  \quad 1\le i \le k.
\end{eqnarray}
Then $\xi_{g_i} := \sum_{m=0}^\infty {\cal P}_0 g_i(X_m) \in L^2$
and
\begin{eqnarray}
\label{eq:finiclt} n^{1/2} (P_n-P){\bf g} = n^{-1/2} \sum_{l=1}^n
\{ {\bf g}(X_l) - \E[ {\bf g}(X_1)]\} \Rightarrow N(0, \Sigma)
\end{eqnarray}
where $\Sigma_{i,j} = \E ( \xi_{g_i} \xi_{g_j}) = \sum_{l \in \Z}
{\rm cov}[g_i(X_0), g_j(X_l)]$.
\end{lemma}

\noindent{\it Proof.} The case in which $k=1$ easily follows from
Hannan (1973, Theorem 1) and Woodroofe (1992). For $k \ge 2$, we
apply the Cram\'er-Wold device. Let $\lambda_1, \lambda_2$ be two
real numbers. Using the relation $\xi_{g_i} = \sum_{m=0}^\infty
{\cal P}_0 g_i(X_m)$, we have $\lambda_1 \xi_{g_1} + \lambda_2
\xi_{g_2} = \xi_{\lambda_1 g_1 + \lambda_2 g_2}$, from which
(\ref{eq:finiclt}) easily follows with the stated covariance
function. \qed

There are many other forms of central limit theorems for dependent
random variables. Lemma \ref{lem:hannan} imposes simple and easily
verifiable conditions. In addition, it also provides a very
natural vehicle for the finite-dimensional convergence of
$\{M_n(g):~ g \in {\cal G}_{\gamma, \mu}\}$; see the proof of
Theorem \ref{th:Pn-Pg} in Section \ref{sec:thm3pf}.

Gordin (1969) obtained a general central limit theorem. However,
as pointed out by Hall and Heyde (1980, p. 130), the condition
imposed in Gordin (1969) is very difficult to check. To overcome
the difficulty, in their book Hall and Heyde proposed Theorem 5.3
(see p. 133), which states that the central limit theorem holds
provided $\sum_{n=0}^\infty {\cal P}_0 g(X_n)$ converges to some
random variable $\xi$ in $L^2$ with $\| \xi \| > 0$ and $\lim_{n
\to \infty} n \| (P_n-P)g \|^2 \to \| \xi \|^2$. Due to the
dependence, the verification of the latter condition could be
difficult as well. In Lemma \ref{lem:hannan} we do not need to
verify the latter condition. Another related result is given in
Dedecker and Rio (2000). A key condition in the latter paper is
that $\sum_{n=0}^\infty g(X_0) \E[ g(X_n) | {\cal F}_0]$ converges
in $L^1$. Consider the linear process $g(X_n) = \sum_{i=0}^\infty
a_i \varepsilon_{n-i}$, where $a_i = (i+1)^{-\beta}$, $\beta >
1/2$ and $\varepsilon_i$ are iid standard normal. Then their $L^1$
convergence condition requires $\beta > 3/2$. In comparison,
(\ref{eq:bernoulli}) only needs $\beta > 1$.

\section{Piecewise differentiable functions}
\label{sec:piecewise} Using the weak convergence theory for step
and absolutely continuous functions, we can easily deal with
piecewise differentiable functions. Let $I \ge 1$ be a fixed
integer and ${\cal G}^I_\gamma = \{ \sum_{i=1}^I g_i(\cdot) {\bf
1}_{\cdot \le \theta_i}: ~ g_i \in {\cal G}_{\gamma, 1},
~\theta_i\in \R\}$.

\begin{theorem}
\label{th:mix} Let $\gamma > 0$. Assume that $\E (|X_1|^{\gamma +
\delta}) < \infty$ and $\sup_u f_\varepsilon(u | {\cal F}_0) \le
\tau$ hold for some $\delta > 0$ and $\tau < \infty$. Further
assume that $\sigma(f_\varepsilon, w_{1 + \gamma}) < \infty$ and
$\sigma(f'_\varepsilon, w_{-1 -\gamma} ) < \infty$. Then $\{\sqrt
n (P_n - P)g: ~ g \in {\cal G}^I_\gamma\}$ converges weakly to a
Gaussian process.
\end{theorem}
\noindent{\it Proof.} Without loss of generality let $I = 1$. For
all $\theta \in \R$, the function $g_\theta(\cdot) = g(\cdot) {\bf
1}_{ \cdot \le \theta} + g(\theta) {\bf 1}_{\cdot > \theta} \in
C_\gamma {\cal G}_{\gamma, 1}$ for some constant $C_\gamma <
\infty$. By Theorem \ref{th:Pn-Pg}, under the proposed conditions,
$\{ \sqrt n (P_n - P)g: g \in C_\gamma {\cal G}_{\gamma, 1}\}$
converges weakly. It then suffices to show that $\{g(\theta) \sqrt
n (P_n - P) {\bf 1}_{\cdot > \theta}:~ g \in {\cal G}_{\gamma, 1},
~\theta\in \R\} $ is tight. Recall $R_n(s) = \sqrt n [F_n(s) -
F(s)]$. Notice that $(P_n - P) {\bf 1}_{\cdot > \theta} = - (P_n -
P) {\bf 1}_{\cdot \le \theta}$. By Lemma \ref{lem:supHweighted},
$g^2(\theta) \le C_\gamma (1 + |\theta| )^\gamma$. Hence
\begin{eqnarray*}
\{g(\theta) \sqrt n (P_n - P) {\bf 1}_{\cdot > \theta}:~ g \in
{\cal G}_{\gamma, 1}, ~\theta\in \R\}
 \subset \{ \lambda R_n(s): ~ \lambda^2 \le C_\gamma
(1 + |s|)^\gamma, ~ s\in \R \}.
\end{eqnarray*}
The latter is a process indexed by both $\lambda$ and $s$. By
(\ref{eq:Gnrtail}) and (\ref{eq:Qnrtail}),
\begin{eqnarray}
\label{eq:mixoutr} \lim_{r \to \infty} \limsup_{n \to \infty} \E
\left\{ \sup_{ |s| \ge r} \sup_{|\lambda| \le C^{1/2}_\gamma (1 +
|s|)^{\gamma /2 } } [\lambda^2 R^2_n(s)] \right \}
 = 0.
\end{eqnarray}
Let $\Gamma_r = \{ (\lambda, s): ~ \lambda^2 \le C_\gamma (1 +
|s|)^\gamma, ~ |s| \le r\}$ and $\Gamma = \Gamma_\infty$. For
$(\lambda_1, s_1), (\lambda_2, s_2) \in \Gamma_r$,
\begin{eqnarray*}
|\lambda_1 R_n(s_1) - \lambda_2 R_n(s_2) |
 &\le& |\lambda_1 - \lambda_2| |R_n(s_1)|
  + |\lambda_2| |R_n(s_1) - R_n(s_2)| \\
 &\le& |\lambda_1 - \lambda_2| \sup_{u \in \R} |R_n(u)|
  + C_{r, \gamma} |R_n(s_1) - R_n(s_2)|.
\end{eqnarray*}
Since $\| \sup_{u \in \R} |R_n(u)| \| = {\cal O}(1)$ and
$R_n(\cdot)$ is tight, it is easily seen that $\{ \lambda R_n(s):
~ (\lambda, s) \in \Gamma_r \}$ is also tight. By
(\ref{eq:mixoutr}), the process $\{ \lambda R_n(s): ~ (\lambda, s)
\in \Gamma \}$ is tight.  Notice that
\begin{eqnarray*}
\|{\cal P}_0 g(X_n) {\bf 1}_{X_n \le \theta} \|
 \le \|{\cal P}_0 g_\theta(X_n)\|
 + |g(\theta)| \|{\cal P}_0 {\bf 1}_{X_n \le \theta} \|
\end{eqnarray*}
is summable, the finite-dimensional convergence follows from Lemma
\ref{lem:hannan}. \qed

\section{Applications}
\label{sec:cond1} Recall (\ref{eq:sigmahm}) for the definition of
$\sigma(h, m)$. To apply Theorems \ref{th:RnW}, \ref{th:modulus}
and \ref{th:Pn-Pg}, one needs to verify the short-range dependence
condition that $\sigma(h, m)$ is finite. In many important
applications including Markov chains and linear processes, there
is an $\sigma({\cal F}_n)$-measurable random variable $Y_n$ such
that
\begin{eqnarray}
\label{eq:FY} \p(X_{n+1} \le x |{\cal F}_n) = \p(X_{n+1} \le x
|Y_n).
\end{eqnarray}
Write $Y_n = I(\ldots, \varepsilon_{n-1}, \varepsilon_n)$ and
$Y_n^* = I(\ldots, \varepsilon_{-1}, \varepsilon^*_0,
\varepsilon_1, \ldots, \varepsilon_n)$, where
$(\varepsilon^*_i)_{i \in \Z}$ is an iid copy of
$(\varepsilon_i)_{i \in \Z}$, and $h(\theta, {\cal F}_n) =
h(\theta, Y_n)$. It turns out that $\sigma(h, m)$ is closely
related to a weighted distance between $Y_n$ and $Y_n^*$. Define
\begin{eqnarray}
H_m(y) = \int_\R \left| {\partial \over {\partial y}} h(\theta, y)
\right|^2 m(d \theta).
\end{eqnarray}

\begin{proposition}\label{prop:rho35}
Let $\rho_m(a, b) = | \int_a^b H^{1/2}_m(y) d y |$. Then for $n
\ge 0$,
\begin{eqnarray}
\label{eq:Vhm}
 \int_\R \| {\cal P}_0 h(\theta, Y_n) \|^2 m(d \theta)
 \le \| \rho_m(Y_n, Y_n^*) \|^2.
\end{eqnarray}
Hence we have $\sigma(h, m) < \infty$ if
\begin{eqnarray}
\label{eq:Nn_srd_rho}
 \sum_{n=0}^\infty \|\rho_m(Y_n, Y_n^*) \| < \infty.
\end{eqnarray}
\end{proposition}

\noindent{\it Proof.} Observe that $ \E[ h(\theta, Y_n) | {\cal
F}_{-1}] = \E[ h(\theta, Y^*_n) | {\cal F}_{-1}] = \E[ h(\theta,
Y^*_n) | {\cal F}_{0}]$. Then we have $ {\cal P}_0 h(\theta, Y_n)
= \E[ h(\theta, Y_n) - h(\theta, Y^*_n) | {\cal F}_{0}]$. By the
Cauchy-Schwarz inequality,
\begin{eqnarray*}
\|{\cal P}_0 h(\theta, Y_n)\|^2 \le \|h(\theta, Y_n) - h(\theta,
Y^*_n)\|^2 \le \E \left[ \int_{Y_n}^{Y_n^*} \left| {\partial \over
{\partial y}} h(\theta, y) \right| d y \right]^2.
\end{eqnarray*}
Let $\lambda(y) = \sqrt{H_m(y)}$. Again by the Cauchy-Schwarz
inequality,
\begin{eqnarray*}
\int_\R \|{\cal P}_0 h(\theta, Y_n)\|^2 m(d \theta)
 &\le& \int_\R \E \left[ \int_{Y_n}^{Y_n^*} {1\over {\lambda(y)}}
 \left| {\partial \over {\partial y}} h(\theta, y)
 \right|^2 d y \times \int_{Y_n}^{Y_n^*} \lambda(y) d y \right]
 m(d \theta)\\
 &=& \E \left[ \int_{Y_n}^{Y_n^*} {1\over {\lambda(y)}}
 \int_\R \left| {\partial \over {\partial y}} h(\theta, y)
 \right|^2  m(d \theta) d y
 \times \int_{Y_n}^{Y_n^*} \lambda(y) d y \right] \\
 &=& \E \left[\int_{Y_n}^{Y_n^*} \lambda(y) d y \right]^2.
\end{eqnarray*}
So (\ref{eq:Vhm}) follows. \qed

Let $h(\theta, Y_n) = f_\varepsilon(\theta | Y_n)$. Then $H_m(y)$
can be interpreted as a measure of "local dependence" of $X_{n+1}$
on $Y_n$ at $y$. If $X_{n+1}$ and $Y_n$ are independent, then
$f_\varepsilon (\theta | y)$ does not depend on $y$, and hence
$H_m = 0$. Let
\begin{eqnarray*}
A(y; \delta) = \int_\R [f_\varepsilon(\theta | y) -
f_\varepsilon(\theta | y + \delta)]^2 m(d \theta).
\end{eqnarray*}
Then $\sqrt {A(y; \delta)}$ is the weighted $L^2$ distance between
the conditional densities of the conditional distributions
$[X_{n+1} | Y_n = y]$ and $[X_{n+1} | Y_n = y + \delta]$. Under
suitable regularity conditions, $\lim_{\delta \to 0} \delta^{-2}
A(y; \delta) = H_m(y)$. Therefore the function $H_m$ quantifies
the dependence of $X_{n+1}$ on $Y_n$. Intuitively,
(\ref{eq:Nn_srd_rho}) suggests short-range dependence in the sense
that if we change $\varepsilon_{0}$ in $Y_n$ to $\varepsilon^*_0$,
then the cumulative impact of the corresponding changes in $Y_n$,
measured by the weighted distance $\rho_m$, is finite.

\begin{remark}{\rm
The random variable $Y_n^*$ can be viewed as a coupled version of
$Y_n$. The coupling method is popular. See Section 3.1 in Dedecker
and Prieur (2003a) and Section 4 in their (2003b) for some recent
contributions. \qed}
\end{remark}

\subsection{Iterated random functions}
\label{sec:ifs} Many nonlinear time series models assume the form
of iterated random functions [Elton (1990), Diaconis and Freedman
(1999)]. Let
\begin{eqnarray}
\label{eq:ifs} X_n = R(X_{n-1}, {\varepsilon_n}),
\end{eqnarray}
where $\varepsilon$, $\varepsilon_k$, $k \in \Z$, are iid random
elements and $R(\cdot, \cdot)$ is a bivariate measurable function.
For the process (\ref{eq:ifs}), due to the Markovian property,
(\ref{eq:FY}) is satisfied with $Y_n = X_n$. The existence of
stationary distribution of (\ref{eq:ifs}) has been widely studied
and there are many versions of sufficient conditions; see Diaconis
and Freedman (1999), Meyn and Tweedie (1994), Steinsaltz (1999),
Jarner and Tweedie (2001), Wu and Shao (2004) among others. Here
we adopt the simple condition proposed by Diaconis and Freedman
(1999). Let the Lipschitz constant
\begin{eqnarray*}
L_\varepsilon = \sup_{x \not = x'} { {|R(x, \varepsilon) - R(x',
\varepsilon)|} \over {|x - x'|}}.
\end{eqnarray*}
Assume that there exist $\alpha > 0$ and $x_0$ such that
\begin{eqnarray}
\label{eq:ifscond1} \E[ L_\varepsilon^\alpha + |x_0 - R(x_0,
\varepsilon)|^\alpha] < \infty \mbox{ and } \E [ \log(
L_\varepsilon)] < 0.
\end{eqnarray}
Then there is a unique stationary distribution (Diaconis and
Freedman, 1999). The latter paper also gives a convergence rate of
an arbitrary initial distribution towards the stationary
distribution. Wu and Woodroofe (2000) pointed out that the simple
sufficient condition (\ref{eq:ifscond1}) also implies the {\it
geometric-moment contraction}: There exist $\beta > 0$, $r\in
(0,1)$ and $C < \infty$ such that
\begin{eqnarray}
\label{eq:gmc} \E[ |F(\ldots, \varepsilon_{-1}, \varepsilon_0,
\varepsilon_{1}, \ldots, \varepsilon_n) - F(\ldots,
\varepsilon'_{-1}, \varepsilon'_0, \varepsilon_{1}, \ldots,
\varepsilon_n)|^\beta ] \le C r^n
\end{eqnarray}
holds for all $n \in \N$, where $(\varepsilon_k')_{k \in \Z}$ is
an iid copy of $(\varepsilon_k)_{k \in \Z}$. Hsing and Wu (2004)
and Wu and Shao (2004) argued that (\ref{eq:gmc}) is a convenient
condition to establish limit theorems. Dedecker and Prieur (2003a,
b) discussed the relationship between (\ref{eq:gmc}) and some new
dependence coefficients. Recently Douc et al (2004) considered
subgeometric rates of convergence of Markov chains. In their paper
they adopted total variational distance, while (\ref{eq:gmc})
involves the Euclidean distance.

Bae and Levental (1995) considered empirical central limit
theorems for stationary Markov chains. It seems that the
conditions imposed in their paper are formidably restrictive.
Consider the Markov chain $X_{n+1} = a X_n + \varepsilon_{n+1}$,
where $0< |a| < 1$ and $\varepsilon_i$ are iid standard normal
random variables. Then the transition density is $q(u | v) = (2
\pi)^{-1/2} \exp[- (u - a v)^2/2]$ and the chain has the
stationary distribution $N[0, (1-a^2)^{-1}]$. Let $\alpha (u) = (2
\pi)^{-1/2} (1- a^2)^{1/2} \exp[- u^2 (1-a^2) /2]$ be the marginal
density. It is easily seen that $\sup_{u, v \in \R} [q(u | v) /
\alpha(u)] = \infty$. Conditions (1.1) and (1.2) in Bae and
Levental (1995) require that this quantity is finite. For this
process, by Theorem \ref{th:lp}, $\{\sqrt n (P_n - P) g: ~ g \in
{\cal G}_{\gamma, \mu}\}$ converges weakly to a tight Gaussian
process for any $\gamma \ge 0$ and $\mu \le 1$.

\subsubsection{AR models with ARCH Errors}
Autoregressive models with conditional heteroscedasticity (ARCH)
have received considerable attention over the last two decades.
For recent work see Berkes and Horv\'ath (2004) and Straumann and
Mikosch (2003), where some statistical inference problems of such
models are considered. Here we consider the model
\begin{eqnarray}
\label{eq:black} X_n = \alpha X_{n-1} + \varepsilon_n \sqrt{a^2 +
b^2 X_{n-1}^2},
\end{eqnarray}
where $(\varepsilon_i)_{i \in \Z}$ are iid random variables. In
the case $b = 0$, then (\ref{eq:black}) is reduced to the
classical AR(1) model. If $b \not = 0$, then the conditional
variance of $X_n$ given $X_{n-1}$ is not a constant and the model
is said to be heteroscedastic. We assume without loss of
generality that $b=1$ since otherwise we can introduce $a' = a/b$
and $\varepsilon'_i = \varepsilon_i/b$. Let $R(x, \varepsilon) =
\alpha x + \varepsilon \sqrt{a^2 + x^2}$. Then $L_\varepsilon \le
\sup_x | \partial R(x, \varepsilon) /
\partial x| \le |\alpha| + |\varepsilon|$. Assume $r = \E
[(|\alpha| + |\varepsilon|)^\beta] < 1$ for some $\beta > 0$. Then
(\ref{eq:ifscond1}) holds and there is a unique stationary
distribution. It is easily seen that (\ref{eq:gmc}) holds with
this $\beta$ and $r$, and $\E(|X_1|^\beta) < \infty$.

We now compute $H_{w_\gamma}(y)$. Denote by $f_\varepsilon$ the
density function of $\varepsilon_i$. Write $u = u(\theta, y) =
(\theta -\alpha y) / {\sqrt {a^2 + y^2}}$. Then the conditional
(transition) density $f_\varepsilon(\theta | y) = (a^2 +
y^2)^{-1/2} f_\varepsilon(u)$ and
\begin{eqnarray}\label{eq:devarch}
{\partial \over {\partial y}} f_\varepsilon(\theta | y)
 = -{ {y f_\varepsilon(u)} \over {(a^2 + y^2)^{3/2}}}
 - { {\alpha f'_\varepsilon(u)} \over {a^2 + y^2}}
 + { {({\theta -\alpha y}) y f'_\varepsilon(u)}
 \over {(a^2 + y^2)^2}}.
\end{eqnarray}
Assume that $\kappa_\gamma := \int_\R f^2_\varepsilon(t)
w_\gamma(d t) + \int_\R |f'_\varepsilon(t)|^2 w_{\gamma + 2}(d t)
< \infty$. Elementary calculations show that
\begin{eqnarray*}
\int_\R f^2_\varepsilon( u(\theta, y) ) w_\gamma(d \theta)
 &=& \int_\R f^2_\varepsilon
 \left({ v \over {\sqrt {a^2 + y^2}}} \right)
 (1 + |v + \alpha y|)^\gamma d v \\
 &\le& C_\gamma (1 + |\alpha y|)^\gamma (a^2 + y^2)^{1/2}
 \kappa_\gamma
 + C_\gamma (a^2 + y^2)^{(1+\gamma)/2} \kappa_\gamma, \cr
\int_\R |f'_\varepsilon( u(\theta, y) )|^2 w_\gamma(d \theta)
 &\le& C_\gamma (1 + |\alpha y|)^\gamma (a^2 + y^2)^{1/2}
 \kappa_\gamma + C_\gamma (a^2 + y^2)^{(1+\gamma)/2}\kappa_\gamma
\end{eqnarray*}
and
\begin{eqnarray*}
\int_\R ({\theta -\alpha y})^2 |f'_\varepsilon(u)|^2 w_\gamma(d
\theta)
 &\le& C_\gamma (1 + |\alpha y|)^\gamma (a^2 + y^2)^{3/2}
 \kappa_\gamma
 + C_\gamma (a^2 + y^2)^{(3 + \gamma)/2} \kappa_\gamma.
\end{eqnarray*}
Combining the preceding three inequalities, we have by
(\ref{eq:devarch}) that
\begin{eqnarray}\label{eq:645}
H_{w_\gamma}(y)
 \le \int_\R \left| {\partial \over {\partial y}}
 f_\varepsilon(\theta | y) \right|^2 w_\gamma(d y)
 \le C_{a, \gamma, \alpha} (1 + |y|)^{-3 + \gamma}.
\end{eqnarray}

\begin{theorem}
\label{th:Ar+Arch} Let $\E [(|\alpha| + |\varepsilon| )^\gamma] <
1$ for some $\gamma > 0$ and
\begin{eqnarray}
\label{eq:Ar+Archf} \int_\R f^2_\varepsilon(t) w_{\gamma + 1}(d t)
+ \int_\R |f'_\varepsilon(t)|^2 w_{\gamma+ 3}(d t) < \infty.
\end{eqnarray}
Then there exists $\phi \in (0, 1)$ such that $\rho_{w_{\gamma +
1}}(Y_n, Y_n^*) = {\cal O}(\phi^n)$ and hence $\{\sqrt n (P_n - P)
g:~ g\in {\cal G}_{\gamma, 1} \}$ converges weakly to a tight
Gaussian process.
\end{theorem}

\noindent{\it Proof.} Let $r = \E [(|\alpha| + |\varepsilon|
)^\gamma ]$. By (\ref{eq:gmc}), $ \E(|Y_n - Y_n^*|^\gamma) \le C
r^n$ for some constant $C$. By (\ref{eq:645}), $H_{w_{\gamma+1}}
(y) \le C (1 + |y|)^{ \gamma -2}$. Let $\lambda = \gamma / 2$. If
$0 < \lambda \le 1$, by Lemma \ref{lem:ineqarch},
\begin{eqnarray*}
\rho_{\gamma + 1}(Y_n, Y_n^*)
 = {\cal O}\left[ \left\|
 \int_{Y_n}^{Y_n^*} (1 + |y|)^{\lambda - 1} d y \right \| \right]
 = {\cal O}[\E(|Y_n - Y_n^*|^{2 \lambda})]^{1/2}
 = {\cal O}( r^{n/2} ).
\end{eqnarray*}
If $\lambda > 1$, let $p = \lambda / (\lambda -1)$. By H\"older's
inequality,
\begin{eqnarray*}
\E[(1 + |Y_n|)^{2 (\lambda -1)} |Y_n - Y_n^*|^2]
 &\le& \{\E[(1 + |Y_n|)^{2 (\lambda -1) p}] \}^{1/p} \times
 \{ \E [|Y_n - Y_n^*|^{2 \lambda}] \}^{1/\lambda} \\
 &=& {\cal O}[(r^n)^{1/\lambda}].
\end{eqnarray*}
Hence
\begin{eqnarray*}
\left\|\int_{Y_n}^{Y_n^*} (1 + |y|)^{\lambda - 1} d y \right\|^2
 &\le& \|[( 1 + |Y_n|)^{\lambda - 1} + (1 + |Y_n^*|)^{\lambda - 1}]
 |Y_n-Y_n^*| \|^2 \\
 &\le&4 \E[(1 + |Y_n|)^{2 (\lambda -1)} |Y_n - Y_n^*|^2]
 = {\cal O}[(r^n)^{1/\lambda}].
\end{eqnarray*}
Let $\phi = \max(r^{1/2}, r^{1/(2\lambda)}) < 1$. Then
$\rho_{\gamma + 1}(Y_n, Y_n^*) = {\cal O}(\phi^n)$, and the
theorem follows from (\ref{eq:Nn_srd_rho}), Proposition
\ref{prop:rho35} and (i) of Theorem \ref{th:Pn-Pg}. \qed

\begin{lemma}
\label{lem:ineqarch} Let $0 < \lambda \le 1$. Then for all $u, v
\in \R$, $|\int_v^u (1 + |y|)^{\lambda - 1} d y| \le
{2^{1-\lambda}} |u-v|^{\lambda} / \lambda$.
\end{lemma}
\noindent{\it Proof}. It suffices to consider two cases (i) $u \ge
v \ge 0$ and (ii) $u \ge 0 \ge v$. For case (i),
\begin{eqnarray*}
\int_v^u (1 + |y|)^{\lambda - 1} d y = {1 \over \lambda}[(1 +
u)^\lambda - (1 + v)^\lambda ] \le {1 \over \lambda}(u-v)^\lambda.
\end{eqnarray*}
For the latter case, let $t = u - v$. Then
\begin{eqnarray*}
\int_v^u (1 + |y|)^{\lambda - 1} d y = {1 \over \lambda}[(1 +
u)^\lambda - 1 + (1 + |v|)^\lambda -1 ] \le {1 \over \lambda}[2(1
+ t/2)^\lambda - 2 ] \le { {2^{1-\lambda}} \over \lambda}
t^\lambda
\end{eqnarray*}
and the lemma follows. \qed

\subsection{Linear processes}
\label{sec:lp} Let $X_t= \sum_{i=0}^\infty a_i \varepsilon_{t-i}$,
where $\varepsilon_k$, $k \in \Z$, are iid random variables with
mean $0$ and finite and positive variance, and the coefficients
$(a_i)_{i \ge 0}$ satisfy $\sum_{i=0}^\infty a_i^2 < \infty$.
Assume without loss of generality that $a_0 = 1$. Let
$F_\varepsilon$ and $f_\varepsilon = F_\varepsilon'$ be the
distribution and density functions of $\varepsilon_k$. Then the
conditional density of $X_{n+1}$ given ${\cal F}_n$ is
$f_\varepsilon(x - Y_n)$ and (\ref{eq:FY}) is satisfied, where
$Y_n = X_{n+1} - \varepsilon_{n+1}$. Limit theorems for short and
long-range dependent linear processes are presented in Sections
\ref{sec:lpsrd} and \ref{sec:lplrd} respectively.

\subsubsection{Short-memory linear processes}
\label{sec:lpsrd}
\begin{proposition}
\label{prop:lp} Let $\gamma \ge 0$. Assume $\E(|\varepsilon_k|^{2
+ \gamma}) < \infty$ and
\begin{eqnarray}
\kappa := \int_\R |f'_\varepsilon(u)|^2 w_\gamma(d u) < \infty.
\end{eqnarray}
Then $\rho_{w_\gamma}(Y_n, Y_n^*) = {\cal O}(|a_{n+1}|)$.
\end{proposition}
\noindent{\it Proof}. Using the elementary inequality $1+ |v+y|
\le (1+|v|)(1+|y|)$, we have
\begin{eqnarray*}
H_{w_\gamma}(y)
 &=& \int_\R |f'_\varepsilon(u-y)|^2 (1+|u|)^\gamma d u \\
 &\le& (1+|y|)^\gamma \int_\R |f'_\varepsilon(v)|^2 w_\gamma(d v)
 = \kappa (1+|y|)^\gamma
\end{eqnarray*}
and
\begin{eqnarray*}
\left|\int_{Y_n}^{Y_n^*} \sqrt{H_{w_\gamma}(y)} d y \right|
 &\le& {\sqrt \kappa} \left|\int_{Y_n}^{Y_n^*}
(1+|y|)^{\gamma / 2} d y \right| \\
 &\le& {\sqrt \kappa} [(1+|Y_n|)^{\gamma / 2} +
(1+|Y_n^*|)^{\gamma / 2}] |Y_n -  Y_n^*|.
\end{eqnarray*}
Observe that ${Y_n}$ and ${Y_n^*}$ are identically distributed.
Then
\begin{eqnarray*}
\rho_{w_\gamma}(Y_n, Y_n^*)
 &\le& 2{\sqrt \kappa} \|(1+|Y_n|)^{\gamma / 2}|Y_n -  Y_n^*| \|\\
 &=& 2{\sqrt \kappa} |a_{n+1}| \|(1+|Y_n|)^{\gamma / 2}
 (\varepsilon_0 - \varepsilon'_0) \| \\
 &\le & 2{\sqrt \kappa} |a_{n+1}| [\|(1+|Y_n|)^{\gamma / 2}
 \varepsilon_0 \|
 + \|(1+|Y_n|)^{\gamma / 2} \varepsilon'_0 \|.
\end{eqnarray*}
Since $\E(|\varepsilon_k|^{2 + \gamma}) < \infty$, it is easily
seen that $\|(1+|Y_n|)^{\gamma / 2} \varepsilon'_0 \| = {\cal
O}(1)$ and $\|(1+|Y_n|)^{\gamma / 2} \varepsilon_0 \| = {\cal
O}(1)$. So the proposition follows. \qed

\begin{remark}{\rm
By Proposition \ref{prop:lp}, if $\gamma = \mu = 0$,
$\varepsilon_k$ has second moment and $\int_\R [ f'_\varepsilon(t)
]^2 d t < \infty$, then (\ref{eq:srdG}) is equivalent to
\begin{eqnarray}
\label{eq:srdlp} \sum_{n = 1}^\infty |a_n| < \infty,
\end{eqnarray}
which is a well-known condition for a linear process to be
short-range dependent.
}\end{remark}

By Theorem \ref{th:Pn-Pg}, Propositions \ref{prop:rho35} and
\ref{prop:lp}, we have
\begin{theorem}
\label{th:lp} Let $\gamma \ge 0$ and $0 \le \mu \le 1$. Assume
(\ref{eq:srdlp}), $\E(|\varepsilon_1|^{2 + \gamma + \mu}) <
\infty$ and
\begin{eqnarray}
\int_\R |f'_\varepsilon(u)|^2 w_{\gamma + \mu}(d u) < \infty.
\end{eqnarray}
Then $\{\sqrt n (P_n - P) g: ~ g \in {\cal G}_{\gamma, \mu}\}$
converges weakly to a tight Gaussian process.
\end{theorem}

\begin{corollary}
\label{cor:EmpWtdlp} Let $\gamma \ge 0$. Assume (\ref{eq:srdlp}),
$\E(|\varepsilon_1|^{2 + \gamma}) < \infty$, $\sup_u
f_\varepsilon(u) < \infty$ and
\begin{eqnarray*}
\int_\R |f'_\varepsilon(u)|^2 w_{\gamma + 1}(d u)
 + \int_\R |f''_\varepsilon(u)|^2 d u < \infty.
\end{eqnarray*}
Then $\{ R_n(s) ( 1 + |s|)^{ \gamma/2 }, ~ s \in \R \}$ converges
weakly to a tight Gaussian process.
\end{corollary}

The corollary easily follows from Propositions \ref{prop:rho35},
\ref{prop:lp} and Corollary \ref{cor:EmpWeighted}. We omit the
details of the proof.

\subsubsection{Long-memory linear processes}
\label{sec:lplrd} Let $a_0 = 1$ and $a_n = n^{-\beta} \ell(n)$, $n
\ge 1$, where $1/2 < \beta < 1$ and $\ell$ is a slowly varying
function (Feller, 1971, p. 275). Then the covariances are not
summable and we say that $(X_t)_{t \in \Z}$ is {\it long-range
dependent} or {\it long-memory}. Let $K$ be a measurable function
such that $K_\infty(x) := \E[K(X_1 + x)]$ is in ${\cal C}^p$, the
class of functions having up to $p$th order derivatives; let
$\sigma_{n, p}^2= n^{2 - p(2 \beta -1)}\ell^{2 p}(n)$ and
\begin{eqnarray}
\label{eq:SnKp} S_n(K; p) = \sum_{i=1}^n \left [K(X_i) -
\sum_{j=0}^p K_\infty^{(j)}(0) U_{i,j}\right], {\rm\, where\,
}\,\, U_{n,r} = \sum_{0 \le j_1 < \ldots <j_r} \prod_{s=1}^r
a_{j_s} \varepsilon_{n - j_s}.
\end{eqnarray}
Ho and Hsing (1996, 1997) initiated the study of such expansions,
which are closely related to chaotic representations, Volterra
processes (Doukhan, 2003) and Hermite expansions if $X_n$ is
Gaussian. Here we consider the weak convergence of $\{S_n(K; p): ~
K \in {\cal K}_\gamma\}$ for the class ${\cal K}_\gamma = {\cal
G}_{\gamma, 0, 0}^*$; recall (\ref{eq:Gstar}) for the definition
of ${\cal G}_{\gamma, \mu, 0}^*$. If $(1 + p) (2\beta - 1) < 1$,
then the limiting distribution of $\{ \sigma^{ -1}_{n, p+1} S_n(K;
p): ~ K \in {\cal K}_\gamma\}$ is degenerate since it forms a line
of multiples of the multiple Wiener-Ito integral $Z_{p+1, \beta}$
(cf. Corollary 3 in Wu (2003a)). See Major (1981) for the
definition of multiple Wiener-Ito integrals. To have a complete
characterization, one needs to consider the case in which $(1 + p)
(2\beta - 1) > 1$. It turns out that, with the help of Proposition
\ref{prop:Mngtight}, we are able to show that the limiting
distribution is a non-degenerate Gaussian process with a $\sqrt
n$-normalization.

Recall $Y_{n-1} = X_n - \varepsilon_n$. Let $S^0_n(y; p) =
\sum_{i=1}^n L^0_p({\cal F}_i, y)$, where
\begin{eqnarray}
\label{eq:HH96L0} L^0_p({\cal F}_n, y) = F_\varepsilon(y -
Y_{n-1}) - F(y) - \sum_{i=1}^p (-1)^i F^{(i)}(y) U_{n,i}.
\end{eqnarray}
Let $L_p({\cal F}_n, y) = [{\bf 1}_{X_n \le y} - F_\varepsilon(y -
Y_{n-1})] + L^0_p({\cal F}_n, y)$, $l^0_p({\cal F}_n, y) =
\partial L^0_p({\cal F}_n, y) / \partial y $, $s^0_n(y; p) =
\sum_{i=1}^n l^0_p({\cal F}_i, y)$ and $S_n(y; p) = \sum_{i=1}^n
L_p({\cal F}_i, y)$.

\begin{theorem}
\label{th:pbless1} Assume $\E(| \varepsilon_1 |^{4 + \gamma} ) <
\infty$ for some $\gamma \ge 0$, $f_\varepsilon \in {\cal C}^{
p+1}$ and
\begin{equation}
\label{eq:fkappa1} \sum_{r=0}^{p+1} \int_{\R} |
f_\varepsilon^{(r)}(x)|^2 w_{\gamma}(dx) < \infty.
\end{equation}
(i) If  $p < (2\beta - 1)^{-1} - 1$, then
\begin{eqnarray}\label{eq:MWI}
{1\over {\sigma_{n, p+1}}} \{ S_n(K; p): ~ K \in {\cal K}_\gamma\}
\Rightarrow \{K_\infty^{(p+1)}(0): ~ K \in {\cal K}_\gamma\}
Z_{p+1, \beta}.
\end{eqnarray}
(ii) If  $p > (2\beta - 1)^{-1} -1$, then $\{n^{-1/2} S_n(K; p): ~
K \in {\cal K}_\gamma\}$ converges weakly to a tight Gaussian
process.
\end{theorem}
\noindent{\it Proof}. As in Wu (2003a), let $\theta_{n, p} =
|a_{n-1}| [|a_{n-1}| + A^{1/2}_n(4) + A_n^{p/2}(2)]$ and
$\Theta_{n, p}=\sum_{k=1}^n \theta_{k, p}$, where $A_n(k) =
\sum_{i=n}^\infty |a_i|^k$. By Karamata's theorem (Feller, 1971,
p. 281), $A_n(k) = {\cal O}(n |a_n|^k)$, $k \ge 2$. (i) It follows
from Theorem 1 and Corollary 3 in Wu (2003a). (ii) Since $(1 + p)
(2\beta - 1) > 1$,
\begin{eqnarray}
\label{eq:sumtheta} \sum_{n=1}^\infty \theta_{n, p} =
\sum_{n=1}^\infty {\cal O} [a_n^2 + |a_n| (n a_n^4)^{1/2} + |a_n|
(n a_n^2)^{p/2}] < \infty.
\end{eqnarray}
By Lemma 9 in Wu (2003a), the condition $\sum_{r=0}^{p} \int_{\R}
| f_\varepsilon ^{(r)}(x) |^2 w_{\gamma}(d x) < \infty$ together
with $\E(| \varepsilon_1 |^{4 + \gamma} ) < \infty$ implies that
\begin{eqnarray}
\label{eq:P1Ln} \int_{\R} \|{\cal P}_1 L ({\cal F}_n,t) \|^2
w_{\gamma}(d t) = {\cal O} (\theta_{n, p}^2).
\end{eqnarray}
Using the same argument therein, it can be shown that
(\ref{eq:fkappa1}) entails the similar result
\begin{eqnarray}
\label{eq:P1l0n} \int_{\R} \|{\cal P}_1 l^0 ({\cal F}_n, t) \|^2
w_{\gamma}(d t) = {\cal O} (\theta_{n, p}^2).
\end{eqnarray}
By Lemma 8 in Wu (2003a), $\E(| \varepsilon_1 |^{\max(1 + \gamma,
2)} ) < \infty$ and (\ref{eq:fkappa1}) imply that $S_n(K;p)$ has
the representation
\begin{eqnarray}
\label{eq:SnKp2} S_n(K; p) = -\int_\R K'(x) S_n(x; p) d x = -
\sum_{i=1}^n \int_\R K'(x) L_p({\cal F}_i, x ) d x.
\end{eqnarray}
Combining (\ref{eq:P1Ln}) and (\ref{eq:sumtheta}), since $K \in
{\cal K}_\gamma$,
\begin{eqnarray*}
\| {\cal P}_0 \int_\R K'(y) L ({\cal F}_n, y) d y \|
 &\le& \sqrt{\int_\R [K'(y)]^2 w_{-\gamma}(d y)
 \int_{\R} \|{\cal P}_0 L ({\cal F}_n, t) \|^2 w_{\gamma}(d t)}
 \\
 &=& {\cal O}(\theta_{n+1, p})
\end{eqnarray*}
are summable and by Lemma \ref{lem:hannan} the finite-dimensional
convergence follows.

We now use the truncation argument as in the proofs of
Propositions \ref{prop:Mngtight} and \ref{prop:Nngtight} to
establish the tightness. Since ${\cal P}_0 L ({\cal F}_n, t) =
{\cal P}_0 L^0 ({\cal F}_n, t)$ for $n \not = 0$, (\ref{eq:P1Ln})
implies that $\int_{\R} \|{\cal P}_1 L^0 ({\cal F}_n,t) \|^2
w_{\gamma}(d t) = {\cal O} (\theta_{n, p}^2)$. By Lemma
\ref{lem:Intsum} and (\ref{eq:sumtheta}), we have
\begin{eqnarray}
\label{eq:S0} \lim_{r \to \infty} \limsup_{n \to \infty} {1 \over
n} \int_{ |x| > r} \| S^0_n(x; p) \|^2 w_{\gamma}(d x) = 0,
\end{eqnarray}
and
\begin{eqnarray}
\label{eq:s0} \limsup_{n \to \infty} {1 \over n} \int_{ \R } \|
s^0_n(x; p) \|^2 w_{\gamma}(d x) < \infty.
\end{eqnarray}
Applying the maximal inequality (\ref{eq:SupHweighted}) of Lemma
\ref{lem:supHweighted} with $\mu = 0$,
\begin{eqnarray}
\label{eq:supS0} {1 \over n} \E\left[\sup_x |S^0_n(x; p)|^2 (1 +
|x|)^{\gamma}\right] \le {1 \over n} \int_{ \R } [ \| S^0_n(x; p)
\|^2 +  \| s^0_n(x; p) \|^2] w_{\gamma}(d x) = {\cal O}(1).
\end{eqnarray}
Note that $S_n(x; p) = n[F_n(x) - {\tilde F}_n(x)] + S^0_n(x; p)$.
By (\ref{eq:SnKp2}),
\begin{eqnarray*}
S_n(K; p) &=& -\int_\R K'(x) n[F_n(x) - F(x)] d x
     -\int_\R K'(x) S^0_n(x; p) d x \cr
 &=& \sqrt n M_n(K) - \int_{ |x|>r } K'(x) S^0_n(x; p) d x
 -\int_{-r}^r K'(x) S^0_n(x; p) d x.
\end{eqnarray*}
Since $\E( |X_1|^{\gamma + 1}) < \infty$, by (i) of Proposition
\ref{prop:Mngtight}, the process $\{M_n(K): ~ K \in {\cal
K}_\gamma\}$ is tight. For the second term, by the Cauchy-Schwarz
inequality,
\begin{eqnarray*}
{1\over n} \left\| \int_{ |x| > r} K'(x) S^0_n(x; p) d x
\right\|^2 \le \int_{ |x| > r} [K'(x)]^2 w_{-\gamma}(d x) \times
{1\over n} \int_{ |x| > r} \| S^0_n(x; p) \|^2 w_{\gamma}(d x),
\end{eqnarray*}
which converges to $0$ by first letting $n \to \infty$ and then
$r\to \infty$. Now we deal with the last term. Using integration
by parts,
\begin{eqnarray*}
-\int_{-r}^r K'(x) S^0_n(x; p) d x = \int_{-r}^r K(x) s^0_n(x; p)
d x - K(r) S^0_n(r; p) + K(-r) S^0_n(-r; p)
\end{eqnarray*}
As in the proof of Proposition \ref{prop:Nngtight}, (\ref{eq:s0})
implies that $\{ n^{-1/2} \int_{-r}^r K(x) s^0_n(x; p) d x, ~ K
\in {\cal K}_\gamma\}$ is tight. By (\ref{eq:supS0}), $\|S^0_n(\pm
r; p)\| = {\cal O}(\sqrt n)$. Since $|K(x)|^2 \le C_{\gamma} ( 1 +
|x|)^{\gamma + 1}$, it is easily seen that $\{ n^{-1/2} K(\pm r)
S^0_n(\pm r; p), ~ K \in {\cal K}_\gamma\}$ is also tight. \qed

\bigskip

It is unclear how to generalize Theorem \ref{th:pbless1} to
long-range dependent heavy-tailed processes, linear fields
(non-causal) and other long-range dependent processes. The special
case $p=1$ is considered in many earlier results; see Giraitis and
Surgailis (1999), Doukhan et al (2002) and Doukhan et al (2004).
Limit theorems for heavy-tailed processes are presented in Hsing
(1999), Koul and Surgailis (2001), Surgailis (2002), Wu (2003b)
and Pipiras and Taqqu (2003). It seems that there is no easy way
to establish $p$th order expansions for $p \ge 2$.

\section{Inequalities}
\label{sec:inequalities} The inequalities presented in this
section are of independent interest and they may have wider
applicability. They are used in the proofs of the results in other
sections.

\begin{lemma}
\label{lem:supHweighted} Let $H \in {\cal AC}$. (i) Let $\mu \le
1$ and $\gamma \in \R$. Then there exists $C_{\gamma, \mu} <
\infty$ such that
\begin{eqnarray}
\label{eq:SupHweighted}
 \sup_{x \in \R} [H^2(x) (1+|x|)^{\gamma}]
 \le C_{\gamma, \mu} \int_\R H^2(u) w_{\gamma - \mu} (d u)
 + C_{\gamma, \mu} \int_\R [H'(u)]^2 w_{\gamma + \mu} (d u).
\end{eqnarray}
(ii) Let $\gamma > 0$, $\mu = 1$ and $H(0) = 0$. Then
\begin{eqnarray}
\label{eq:SupHmu1}
 \sup_{x \in \R} [H^2(x) (1+|x|)^{-\gamma}]
 \le {1 \over \gamma} \int_\R [H'(u)]^2 w_{-\gamma + 1} (d u)
\end{eqnarray}
and
\begin{eqnarray}
\label{eq:IntHmu1} \int_\R H^2(u) w_{- \gamma - 1} (d u)
 \le {4 \over {\gamma^2}} \int_\R [H'(u)]^2 w_{- \gamma + 1} (d u).
\end{eqnarray}
(iii) Let $\gamma > 0$ and $H(\pm \infty) = 0$. Then
\begin{eqnarray}
\label{eq:SupHmu1infty}
 \sup_{x \in \R} [H^2(x) (1+|x|)^{\gamma}]
 \le {1 \over \gamma} \int_\R [H'(u)]^2 w_{\gamma + 1} (d u)
\end{eqnarray}
and
\begin{eqnarray}
\label{eq:IntHmu1infty} \int_\R H^2(u) w_{\gamma - 1} (d u)
 \le {4 \over {\gamma^2}} \int_\R [H'(u)]^2 w_{\gamma + 1} (d u).
\end{eqnarray}
\end{lemma}
\noindent{\it Proof.} (i) By Lemma 4 in Wu (2003a), for $t \in \R$
and $\delta > 0$ we have
\begin{eqnarray}
\label{eq:Lem4Bernoulli}
 \sup_{t\le s \le t+ \delta} H^2(s) \le {2\over \delta}
 \int_t^{t+\delta} H^2(u) d u
 + 2\delta \int_t^{t+\delta} [H'(u)]^2 d u.
\end{eqnarray}
We first consider the case $\mu < 1$. Let $\alpha = 1/(1-\mu)$. In
(\ref{eq:Lem4Bernoulli}) let $t = t_n = n^\alpha$ and $\delta =
\delta_n = (n + 1)^\alpha - n^\alpha$, $n \in \N$ and $I_n = [t_n,
t_{n+1}]$. Since $\delta_n \sim \alpha n^{\alpha - 1}$ as $n \to
\infty$,
\begin{eqnarray}
\label{eq:SupH55} \sup_{x\in I_n} [H^2(x) (1+x)^{\gamma}]
 &\le& 2  \sup_{x\in I_n} (1+x)^{\gamma}
 \left[\delta_n^{-1} \int_{I_n}
 H^2(u) d u +  \delta_n \int_{I_n} [H'(u)]^2 d u \right] \cr
 &\le& C \int_{I_n}
 H^2(u) w_{\gamma - \mu} (d u) + C \int_{I_n}
 [H'(u)]^2 w_{\gamma + \mu} (d u).
\end{eqnarray}
It is easily seen in view of (\ref{eq:Lem4Bernoulli}) that
(\ref{eq:SupH55}) also holds for $n = 0$ by choosing a suitable
$C$. By summing (\ref{eq:SupH55}) over $n=0, 1, \ldots$, we obtain
(\ref{eq:SupHweighted}) with $\sup_{x \in \R}$ replaced by
$\sup_{x \ge 0}$. The other side $x \le 0$ similarly follows.

If $\mu = 1$, we let $t_n = 2^n$, $\delta_n = t_{n+1} - t_n = t_n$
and $I_n = [t_n, t_{n+1}]$, $n = 0, 1, \ldots$. The argument above
similarly yields the desired inequality.

(ii) Let $s \ge 0$. Since $H(s) = \int_0^s H'(u) d u$, by the
Cauchy-Schwarz inequality, (\ref{eq:SupHmu1}) follows from
\begin{eqnarray*}
H^2(s) &\le& \int_0^s |H'(u)|^2 (1 + u)^{1-\gamma} d u \times
\int_0^s (1 + u)^{\gamma -1} d u \\
 &\le& \int_\R [H'(u)]^2
 w_{-\gamma + 1} (d u) \times { { (1+ s)^\gamma -1} \over \gamma}.
\end{eqnarray*}
Applying Theorem 1.14 in Opic and Kufner (1990, p. 13) with $p = q
= 2$, the Hardy-type inequality (\ref{eq:IntHmu1}) easily follows.
The proof of (iii) is similar as (ii).\qed

\bigskip

\begin{lemma}
\label{lem:Intsum} Let $m$ be a measure on $\R$, $A\subset \R$ a
measurable set and $T_n(\theta) = \sum_{i=1}^n h(\theta, {\cal
F}_i)$, where $h$ is a measurable function. Then
\begin{eqnarray}
\label{eq:Intsum} \sqrt{\int_A \|T_n(\theta) - \E [T_n(\theta)]
\|^2 m(d \theta) }
 \le {\sqrt n} \sum_{j=0}^\infty \sqrt{ \int_A
 \| {\cal P}_0 h(\theta, {\cal F}_j) \|^2 m(d \theta) }.
\end{eqnarray}
\end{lemma}
\noindent{\it Proof.} For $j=0, 1, \ldots$ let $T_{n, j}(\theta) =
\sum_{i=1}^n \E[h(\theta, {\cal F}_i) | {\cal F}_{i-j}]$ and $
\lambda^2_j = \int_A \| {\cal P}_0 h(\theta, {\cal F}_j) \|^2 m(d
\theta)$, $\lambda_j \ge 0$. By the orthogonality of $
\E[h(\theta, {\cal F}_i) | {\cal F}_{i-j}] - \E[h(\theta, {\cal
F}_i) | {\cal F}_{i-j-1}]$, $i=1, 2, \ldots, n$,
\begin{eqnarray*}
\int_A \|T_{n, j}(\theta) - T_{n, j-1}(\theta)\|^2 m(d \theta)
 &=& n \int_A \| \E[h(\theta, {\cal F}_1) | {\cal F}_{1-j}] -
 \E[h(\theta, {\cal F}_1) | {\cal F}_{-j}] \|^2 m(d \theta) \\
 &=& n \int_A \| {\cal P}_{1-j} h(\theta, {\cal F}_1) \|^2
 m(d \theta) = n \lambda_j^2.
\end{eqnarray*}
Note that $T_n(\theta) = T_{n, 0}(\theta)$. Let $\Delta =
\sum_{j=0}^\infty \lambda_j$. By the Cauchy-Schwarz inequality,
\begin{eqnarray*}
\int_A \E |T_n(\theta) - \E [T_n(\theta)]|^2 m(d \theta)
 &=& \int_A \E  \left\{ \sum_{j=0}^\infty
 [T_{n, j}(\theta) - T_{n, j+1}(\theta)] \right\}^2 m(d \theta) \\
 &\le& \Delta \int_A \E \left\{ \sum_{j=0}^\infty \lambda_j^{-1}
 [T_{n, j}(\theta) - T_{n, j+1}(\theta)]^2 \right\} m(d \theta)
 = n \Delta^2
\end{eqnarray*}
and (\ref{eq:Intsum}) follows. \qed

\bigskip

Lemma \ref{lem:mdsumLq} easily follows from Burkholder's
inequality. We omit the proof.
\begin{lemma}
\label{lem:mdsumLq} Let $(D_i)_{1\le i\le n}$ be $L^q$ ($q > 1$)
martingale differences. Then
\begin{eqnarray}
\label{eq:mdsumLq} \| D_1 + \ldots + D_n \|_q^{\min(q, 2)}
 \le [18 q^{3/2} (q-1)^{-1/2}]^{\min(q, 2)}
 \sum_{i=1}^n \| D_i \|_q^{\min(q, 2)}.
\end{eqnarray}
\end{lemma}

\bigskip

Lemma \ref{lem:maximal} below gives a simple maximal inequality.
Weaker and special versions of it can be found in Doob (1953), Wu
and Woodroofe (2004) and Billingsley (1968). It has the advantage
that the dependence structure of $\{Z_i\}$ can be arbitrary.

\begin{lemma}
\label{lem:maximal} Let $q > 1$ and $Z_i$, $1\le i\le 2^d$, be
random variables in ${\cal L}^q$, where $d$ is a positive integer.
Let $S_n = Z_1 + \ldots + Z_n$ and $S_n^* = \max_{i \le n} |S_i|$.
Then
\begin{eqnarray}
\label{eq:maximal} \| S_{2^d}^* \|_q \le \sum_{r=0}^d \left[
\sum_{m=1}^{2^{d-r}} \| S_{2^r m} - S_{2^r (m-1)} \|_q^q
\right]^{1\over q}.
\end{eqnarray}
\end{lemma}

\noindent{\it Proof}. Let $p = q/(q-1)$ and $\Lambda =
\sum_{r=0}^d \lambda_r^{-p}$, where
\begin{eqnarray*}
\lambda_r = \left[ \sum_{m=1}^{2^{d-r}} \| S_{2^r m} - S_{2^r
(m-1)} \|_q^q \right]^{-1\over {p+q} }.
\end{eqnarray*}
For the positive integer $k \le 2^d$, write its dyadic expansion
$k = 2^{r_1} + \ldots +2^{r_j}$, where $0 \le r_j < \ldots < r_1
\le d$, and $k(i) = 2^{r_1} + \ldots +2^{r_i}$. By H\"older's
inequality,
\begin{eqnarray*}
|S_k|^q &\le& \left[ \sum_{i=1}^j |S_{k(i)} - S_{k(i-1)}| \right]
^q \\
 &\le& \left[ \sum_{i=1}^j \lambda_{r_i}^{-p} \right]^{q/p} \left[
 \sum_{i=1}^j \lambda_{r_i}^{q} |S_{k(i)} - S_{k(i-1)}|^q
 \right] \\
 &\le& \Lambda^{q/p} \sum_{i=1}^j \lambda_{r_i}^{q}
 \sum_{m=1}^{2^{d-r_i}} | S_{2^{r_i} m} - S_{2^{r_i} (m-1)}|^q \\
 &\le& \Lambda^{q/p} \sum_{r=0}^d \lambda_{r}^{q}
 \sum_{m=1}^{2^{d-r}} | S_{2^{r} m} - S_{2^{r} (m-1)}|^q,
\end{eqnarray*}
which entails $\| S_{2^d}^* \|_q^q \le \Lambda^{q/p} \sum_{r=0}^d
\lambda_{r}^{q} \lambda_{r}^{-p-q} = \Lambda^q$ and hence
(\ref{eq:maximal}). \qed.

\section{Proofs of Theorems \ref{th:RnW} and \ref{th:modulus}}
\label{sec:thm12pr} Following (\ref{eq:DecompF}), let $G_n(s) =
n^{1/2} [F_n(x) - {\tilde F}_n(x)]$ and $Q_n(s) = n^{1/2} [{\tilde
F}_n(x) - F(x)]$. Then $R_n(s) = G_n(s) + Q_n(s)$. Sections
\ref{sec:Gn} and \ref{sec:Qn} deal with $G_n$ and $Q_n$
respectively. Theorems \ref{th:RnW} and \ref{th:modulus} are
proved in Sections \ref{sec:83} and \ref{sec:84}.

\subsection{Analysis of $G_n$}
\label{sec:Gn} The main result is this section is Lemma
\ref{lem:Gnr} which concerns the weak convergence of $G_n$.
\begin{lemma}
\label{lem:Bernoulli48} Let $q \ge 2$. Then there is a constant
$C_q < \infty$ such that
\begin{eqnarray}
\label{eq:Bernoulli48}
 \| G_n(y) - G_n(x) \|_q^q
  &\le & C_q n^{\max(1, q/4) - q/2} [F(y) - F(x)] \cr
   & &  + C_q (y-x)^{ q/2 -1}
\int_x^y \E [f_\varepsilon^{q/2} (u | {\cal F}_0)] d u
\end{eqnarray}
holds for all $n \in\N$ and all $x < y$, and
\begin{eqnarray}
\label{eq:Bernoulli49}
 \|G_n(x) \|_q^q \le C_q \min[F(x), 1 - F(x)].
\end{eqnarray}
\end{lemma}

\noindent{\it Proof}.  Let $d_i(s) = {\bf 1}_{X_i \le s} - \E(
{\bf 1}_{X_i \le s} | {\cal F}_{i-1})$, $d_i = d_i(y) - d_i(x)$
and $D_i = d_i^2 - \E( d_i^2 | {\cal F}_{i-1})$. Wu (2003a) deals
with the special case in which $X_i$ is a linear process and $q=4$
[cf. Inequality (48) therein]. Let $q' = q/2$. By Burkholder's
inequality (Chow and Teicher, 1978),
\begin{eqnarray}
\label{eq:Bernoulli48y-x} \| G_n(y) - G_n(x) \|_q^q
 &=& n^{-q/2} \E ( |d_1 + \ldots + d_n|^q) \cr
 &\le& C_q n^{-q/2} \E [ (d_1^2 + \ldots + d_n^2)^{q/2}] \cr
 &\le& C_q n^{-q/2}
   \left\|\sum_{i=1}^n D_i\right\|^{q'}_{q'}
   + C_q n^{-q/2} \E\left\{ \left[ \sum_{i=1}^n
   \E(d_i^2 | {\cal F}_{i-1}) \right]^{q'} \right\}.
\end{eqnarray}
Since $D_i$, $i\in \Z$, form stationary martingale differences, by
Lemma \ref{lem:mdsumLq} we have
\begin{eqnarray}
\label{eq:SumDi} \left\|\sum_{i=1}^n D_i\right\|_{q'}^{q'}
 &\le& C_q n^{\max(1, q'/2)} \|D_1\|_{q'}^{q'} \cr
 &\le& C_q n^{\max(1, q'/2)} 2^{q' - 1} [ \|d_1^2\|_{q'}^{q'}
  + \| \E( d_1^2 | {\cal F}_0) \|_{q'}^{q'} ] \cr
 &\le& C_q n^{\max(1, q'/2)} 2^{q'} \|d_1^2\|_{q'}^{q'},
\end{eqnarray}
where we have applied Jensen's inequality $\| \E( d_1^2 | {\cal
F}_0) \|_{q'}^{q'} \le \| d_1^2 \|_{q'}^{q'}$. Notice that $|d_1|
\le 1$,
\begin{eqnarray}
\label{eq:TermDi} \|d_1^2\|_{q'}^{q'}
 \le \|d_1\|_{q'}^{q'}
 \le 2^{q'-1} [ \|{\bf 1}_{ x\le X_i\le y} \|_{q'}^{q'}
 + \| \E({\bf 1}_{ x\le X_i\le y} | {\cal F}_0) \|_{q'}^{q'}]
 \le 2^{q'} [F(y) - F(x)].
\end{eqnarray}
On the other hand, since $q' \ge 1$ and $\E(d_1^2 | {\cal F}_0)
\le \E({\bf 1}_{ x\le X_i\le y} | {\cal F}_0)$, we have by
H\"older's inequality with $p' = q'/(q'-1)$ that
\begin{eqnarray}
\label{eq:SecondTermDi} \E\left\{ \left[ \sum_{i=1}^n
   \E(d_i^2 | {\cal F}_{i-1}) \right]^{q'} \right\}
 &\le& n^{q'} \| \E(d_1^2 | {\cal F}_0) \|_{q'}^{q'} \cr
 &\le& n^{q'} \E \left\{ \left[ \int_x^y f_\varepsilon ( u | {\cal
 F}_0) d u \right]^{q'} \right\} \cr
 &\le& n^{q'} \E \left[ (y-x)^{q'/p'} \int_x^y
 f_\varepsilon^{q'} ( u | {\cal F}_0) d u \right].
\end{eqnarray}
Combining (\ref{eq:Bernoulli48y-x}), (\ref{eq:SumDi}),
(\ref{eq:TermDi}) and (\ref{eq:SecondTermDi}), we have
(\ref{eq:Bernoulli48}).

To show (\ref{eq:Bernoulli49}), we let in
(\ref{eq:Bernoulli48y-x}) that $d_i = d_i(x) = {\bf 1}_{X_i \le x}
- \E( {\bf 1}_{X_i \le x} | {\cal F}_{i-1})$. Then
\begin{eqnarray*}
\E ( |d_1 + \ldots + d_n|^q)
 &\le& C_q n^{\max(1, q/2)} \|d_1\|_q^q \\
 &\le& C_q n^{q/2} \|d_1\|^2 \le C_q n^{q/2} F(x) [ 1- F(x)]
\end{eqnarray*}
completes the proof. \qed

\bigskip

\begin{lemma}
\label{lem:ossi-emp} Let $q > 2$. Then there exists a constant
$C_q < \infty$ such that
\begin{eqnarray}
\label{eq:ossi-emp}
 \lefteqn{ \E \left[\sup_{0\le s< b} |G_n(a+s)
 - G_n(a)|^q \right]
 \le C_q d^q n^{\max(1, q/4) - q/2} [ F(a+b) - F(a)]} \cr
 & + & C_q b^{q/2 - 1} [1 + n^{q/2} 2^{d (1-q/2)}]
 \int_a^{a+b} \E [f_\varepsilon^{q/2} (u | {\cal F}_0)] d u
\end{eqnarray}
holds for all $b > 0$, $a \in \R$ and $n, d \in \N$. In
particular, for $d = 1+ \lfloor (\log n) / [ (1-2/q) \log 2]
\rfloor$, we have
\begin{eqnarray}
\label{eq:ossi-emp-d}
 \E \left[\sup_{0\le s< b} |G_n(a+s) - G_n(a)|^q \right]
 &\le& C_q (\log n)^q n^{\max(1, q/4) - q/2} [ F(a+b) - F(a)]\cr
 &+& C_q b^{q/2 - 1}
 \int_a^{a+b} \E [f_\varepsilon^{q/2} (u | {\cal F}_0)] d u.
\end{eqnarray}
\end{lemma}

\noindent{\it Proof}. Let $h = b 2^{-d}$, $Z_j = G_n(a + j h) -
G_n(a + (j-1) h)$, $j = 1, \ldots, 2^d$ and $S_j = Z_1 + \ldots +
Z_j$.  By Lemma \ref{lem:Bernoulli48},
\begin{eqnarray*}
\| S_{2^r m} - S_{2^r (m-1)} \|_q^q
 &\le& C_q n^{\max(1, q/4) - q/2}
 [F(a + 2^r m h) - F(a + 2^r (m-1) h)] \cr
 & &  + C_q (2^r h)^{ q/2 -1}
 \int_{a + 2^r (m-1) h}^{a + 2^r m h}
 \E [f_\varepsilon^{q/2} (u | {\cal F}_0)] d u
\end{eqnarray*}
Hence
\begin{eqnarray*}
\sum_{m=1}^{2^{d-r}} \| S_{2^r m} - S_{2^r (m-1)} \|_q^q
 &\le& C_q n^{\max(1, q/4) - q/2} [F(a + b) - F(a)] \\
 &&+ C_q (2^r h)^{ q/2 -1}
\int_{a}^{a + b} \E [f_\varepsilon^{q/2} (u | {\cal F}_0)] d u.
\end{eqnarray*}
By Lemma \ref{lem:maximal},
\begin{eqnarray}
\label{eq:ossi-empZ} \| S_{2^d}^* \|_q
 &\le& \sum_{r=0}^d \{ C_q n^{\max(1, q/4) - q/2} [F(a + b) - F(a)]
 \}^{1/q} \cr
 &&+ \sum_{r=0}^d \left \{ C_q (2^r h)^{ q/2 -1}
 \int_{a}^{a + b} \E [f_\varepsilon^{q/2} (u | {\cal F}_0)] d u
 \right\}^{1/q} \cr
 &\le & d \{ C_q n^{\max(1, q/4) - q/2} [F(a + b) - F(a)]
 \}^{1/q}\cr
 && + \left \{ C_q (2^d h)^{ q/2 -1}
\int_{a}^{a + b} \E [f_\varepsilon^{q/2} (u | {\cal F}_0)] d u
\right\}^{1/q}.
\end{eqnarray}
Recall ${\tilde F}_n (x) = n^{-1} \sum_{i=1}^n F_\varepsilon(x |
{\cal F}_{i-1})$. Let $B_j = \sqrt n [ {\tilde F}_n (a + j h) -
{\tilde F}_n (a + (j-1) h) ]$, $j=1, \ldots, 2^d$, and $q' = q/2$.
Since $0 \le F_\varepsilon \le 1$, by H\"older's inequality,
\begin{eqnarray*}
\| B_j \|_q^q
 &=& n^{-q'} \E \left[ \sum_{i=1}^n \int_{a + (j-1)h}^{a + j h}
 f_\varepsilon(u | {\cal F}_{i-1}) d u \right]^q \\
 &\le& n^{-q'} n^{q-1} \sum_{i=1}^n
 \E \left[\int_{a + (j-1)h}^{a + j h}
 f_\varepsilon(u | {\cal F}_{i-1}) d u \right]^q \\
 &\le& n^{-q'} n^{q-1} \sum_{i=1}^n
 \E \left[\int_{a + (j-1)h}^{a + j h}
 f_\varepsilon(u | {\cal F}_{i-1}) d u \right]^{q'} \\
 &\le& n^{q'} h^{q' -1} \int_{a + (j-1)h}^{a + j h}
 \E [f_\varepsilon^{q/2} (u | {\cal F}_0)] d u.
\end{eqnarray*}
Therefore,
\begin{eqnarray}
\label{eq:ossi-empB} \E\left[\max_{j\le 2^d} B^q_j\right]
 \le \E\left[\sum_{j=1}^{2^d} B^q_j\right]
 \le n^{q'} h^{q' -1} \int_{a}^{a + b}
 \E [f_\varepsilon^{q/2} (u | {\cal F}_0)] d u.
\end{eqnarray}
Observe that
\begin{eqnarray*}
G_n(a+ h \lfloor s/h \rfloor ) - \max_{j\le 2^d} B_j \le G_n(a+s)
\le G_n(a+ h \lfloor s/h + 1 \rfloor ) + \max_{j\le 2^d} B_j.
\end{eqnarray*}
Hence (\ref{eq:ossi-emp}) follows from (\ref{eq:ossi-empZ}),
(\ref{eq:ossi-empB}) and
\begin{eqnarray*}
\sup_{0\le s< b} |G_n(a+s) - G_n(a)|
 &\le &
 \sup_{0\le s< b} |G_n(a+ h \lfloor s/h + 1 \rfloor ) - G_n(a) |
 \\
 & &+ \sup_{0\le s< b} |G_n(a+ h \lfloor s/h \rfloor ) - G_n(a) |
 + 2 \max_{j\le 2^d} B_j \\
 &\le & 2 \max_{j\le 2^d} |G_n(a+ h j) - G_n(a) |
 + 2 \max_{j\le 2^d} B_j \\
 &=& 2 S_{2^d}^* + 2 \max_{j\le 2^d} B_j
\end{eqnarray*}
by noticing that $h = 2^{-d} b$.

For $d = 1+ \lfloor (\log n) / [ (1-2/q) \log 2] \rfloor$, we have
$n^{q/2} 2^{d (1-q/2)} \le 1$ and hence (\ref{eq:ossi-emp-d}) is
an easy consequence of (\ref{eq:ossi-emp}). \qed

\bigskip

\begin{lemma}
\label{lem:Gnr} Let $\gamma \ge 0$ and $q > 2$. Assume $\E[
|X_1|^\gamma + \log(1+ |X_1|)] < \infty$ and (\ref{eq:IntfEq}).
Then (i)
\begin{eqnarray}
\label{eq:Gnr} \E \left[\sup_{s\in \R} |G_n(s)|^q (1+|s|)^\gamma
\right] = {\cal O}(1)
\end{eqnarray}
and (ii) the process $\{ G_n(s)(1+|s|)^{\gamma /q}, ~ s\in \R\}$
is tight and it converges to a tight Gaussian process.
\end{lemma}

\begin{remark} {\rm
In Lemma \ref{lem:Gnr}, the logarithm term $\log(1+ |X_1|)$ is not
needed if $\gamma > 0$. \qed}
\end{remark}

\noindent{\it Proof}. (i) Without loss of generality we show that
$\E [\sup_{s\ge 0} |G_n(s)|^q (1+|s|)^\gamma] = {\cal O}(1)$ since
the case when $s < 0$ follows similarly. Let $\alpha_n = (\log
n)^q n^{\max(1, q/4) - q/2}$. By (\ref{eq:IntfEq}) and
(\ref{eq:ossi-emp-d}) of Lemma \ref{lem:ossi-emp} with $a = b =
2^k$,
\begin{eqnarray}
\label{eq:supG1} &&\sum_{k=1}^\infty (1 + 2^k)^\gamma \E
 \left[\sup_{2^k \le s< 2^{k+1}} |G_n(s) - G_n(2^k)|^q \right] \cr
 &&\le C_q \sum_{k=1}^\infty (1 + 2^k)^\gamma
 \alpha_n [ F(2^{k+1}) - F(2^k)] \cr
 &&+ C_q \sum_{k=1}^\infty (1 + 2^k)^\gamma (2^k)^{q/2 -1}
 \int_{2^k}^{2^{k+1}} \E [f_\varepsilon^{q/2} (u | {\cal F}_0)] d
 u\cr
 && \le C_{\gamma, q} \alpha_n
 \int_2^\infty f(u) (1+ u)^\gamma d u
 + C_{\gamma, q} \int_2^\infty (1 + u)^\gamma u^{q/2 -1}
 \E [f_\varepsilon^{q/2} (u | {\cal F}_0)] d u \cr
 &&\le C_{\gamma, q} \alpha_n
 + C_{\gamma, q} = {\cal O}(1).
\end{eqnarray}
Let $A_\gamma(d) = \sum_{k=1}^d (1+2^k)^\gamma$. Then
$A_\gamma(\lfloor \log_2 v\rfloor +1) \le C (1+v)^\gamma$ if
$\gamma > 0$ and $A_\gamma(\lfloor \log_2 v\rfloor +1) \le \log_2
v$ if $\gamma = 0$. By (\ref{eq:Bernoulli49}) of Lemma
\ref{lem:Bernoulli48}, we have
\begin{eqnarray}
\label{eq:supG2} \sum_{k=1}^\infty (1+2^k)^\gamma \|
G_n(2^k)\|_q^q
 &\le& \sum_{k=1}^\infty
 (1+2^k)^\gamma C_q \int_{2^k}^\infty f(v) d v\cr
 &\le& C_q \int_2^\infty f(v) A_\gamma(\lfloor \log_2 v\rfloor +1)
 d v < \infty.
\end{eqnarray}
Thus (\ref{eq:Gnr}) follows from (\ref{eq:supG1}) and
(\ref{eq:supG2}).

(ii) It is easily seen that the argument in (i) entails
\begin{eqnarray} \label{eq:Gnrtail}
\lim_{r\to \infty} \limsup_{n\to \infty} \E \left[\sup_{|s| > r}
|G_n(s)|^q (1+|s|)^\gamma \right] = 0.
\end{eqnarray}
For $\delta \in(0, 1)$ let the interval $I_k = I_k(\delta) = [k
\delta, (k+1) \delta]$. Then
\begin{eqnarray}\label{eq:diffG}
&&\sup_{s, t\in [-r,r], ~0\le s-t \le \delta } |
 G_n(s)(1+|s|)^{\gamma / q} - G_n(t)(1+|t|)^{\gamma / q}|\cr
 &&\le \sup_{s, t\in [-r,r], ~0\le s-t \le \delta} |
 (1+|s|)^{\gamma / q} [G_n(s) - G_n(t)] | \cr
 &&+ \sup_{s, t\in [-r,r], ~0\le s-t \le \delta} |
 G_n(t)[(1+|s|)^{\gamma / q} - (1+|t|)^{\gamma / q}]|\cr
 &&\le (1 + r)^{\gamma / q}
 \sup_{s, t\in [-r,r], ~0\le s-t \le \delta} |G_n(s) - G_n(t)|
  + C_{r, \gamma, q} \delta \sup_{u\in [-r,r]}|G_n(u)|
\end{eqnarray}
By (i), $\|\sup_{u\in \R}|G_n(u)|\| = {\cal O}(1)$. On the other
hand, by Lemma \ref{lem:ossi-emp},
\begin{eqnarray*}
&&\sum_{k = -\lfloor r/\delta\rfloor -1}^{\lfloor r/\delta \rfloor
+1} \p\left[ \sup_{s \in I_k}
|G_n(s) - G_n(k \delta)| > \epsilon \right]\\
 &&\le \epsilon^{-q} \sum_{k = -\lfloor r/\delta\rfloor
 -1}^{\lfloor r/\delta \rfloor +1}
 \left\{ C_q \alpha_n \p(X_1 \in I_k)
 + C_q \delta^{q/2-1}
 \int_{I_k} \E [f_\varepsilon^{q/2} (u | {\cal F}_0)] d u\right\} \\
 &&\le \epsilon^{-q} C_q \alpha_n +
 \epsilon^{-q} C_q \delta^{q/2-1} \int_{\R}
 \E [f_\varepsilon^{q/2} (u | {\cal F}_0)] d u].
\end{eqnarray*}
By (\ref{eq:IntfEq}), $\int_{\R} \E [f_\varepsilon^{q/2} (u |
{\cal F}_0)] d u < \infty$. Hence
\begin{eqnarray*}
\limsup_{n \to \infty} \p\left[ \sup_{s, t\in [-r,r], ~0\le s-t
\le \delta} |G_n(s) - G_n(t)| > 2 \epsilon\right]
 \le \epsilon^{-q} C_q \delta^{q/2-1},
\end{eqnarray*}
which implies the tightness of $\{G_n(s), ~ -r \le s \le r\}$ for
fixed $r$. So (ii) follows from (\ref{eq:Gnrtail}) and
(\ref{eq:diffG}). \qed

\subsection{Analysis of $Q_n$}
\label{sec:Qn} It is relatively easier to handle $Q_n$ since it is
a differentiable function. The Hardy-type inequalities (cf Lemma
\ref{lem:supHweighted}) are applicable.
\begin{lemma}
\label{lem:Qntight} Assume (\ref{eq:fwep}). Then (i)
\begin{eqnarray}
\label{eq:Qnr} \E \left[\sup_{s\in \R} |Q_n(s)|^2 (1+|s|)^{2\gamma
/ q} \right] = {\cal O}(1)
\end{eqnarray}
and (ii) the process $\{ Q_n(s)(1+|s|)^{\gamma / q}, ~ s\in \R\}$
is tight.
\end{lemma}
\noindent{\it Proof}. Let $\gamma' = 2\gamma / q$. (i) By
(\ref{eq:SupHmu1infty}) of Lemma \ref{lem:supHweighted},
\begin{eqnarray*}
\sup_{|s| \ge r} [Q^2_n(s)(1+|s|)^{\gamma'}]
 \le {1\over {\gamma'}}
 \int_{|s| \ge r} [Q'_n(s)]^2 w_{ 1 + \gamma'}( d s).
\end{eqnarray*}
By Lemma \ref{lem:Intsum},
\begin{eqnarray}
\label{eq:Qnrtail} \left\|\sup_{|s| \ge r} |Q_n(s)|
(1+|s|)^{\gamma'/2} \right\|
 \le {1\over {\sqrt{\gamma'}} }
 \sum_{j=0}^\infty \sqrt{ \int_{|s| \ge r} \| {\cal
 P}_0 f_\varepsilon(\theta|{\cal F}_j) \|^2
 w_{\gamma' + 1}(d \theta) }.
\end{eqnarray}
So (\ref{eq:Qnr}) follows by letting $r = 0$ in
(\ref{eq:Qnrtail}).

(ii) The tightness follows from the similar argument as (ii) of
Lemma \ref{lem:Gnr}. Let $0 < \delta <1$. Then
\begin{eqnarray*}
&&\sup_{s, t\in [-r,r], ~0\le s-t \le \delta } |
 Q_n(s)(1+|s|)^{\gamma / q} - Q_n(t)(1+|t|)^{\gamma / q}|\\
 &&\le \sup_{s, t\in [-r,r], ~0\le s-t \le \delta} |
 (1+|s|)^{\gamma / q} [Q_n(s) - Q_n(t)] | \\
 &&+ \sup_{s, t\in [-r,r], ~0\le s-t \le \delta} |
 Q_n(t)[(1+|s|)^{\gamma / q} -
 (1+|t|)^{\gamma / q}]|\\
 &&\le C_{r, \gamma, q} \delta \sup_{u\in [-r,r]}|Q'_n(u)|
  + C_{r, \gamma, q} \delta \sup_{u\in [-r,r]}|Q_n(u)|
\end{eqnarray*}
Notice that $\|\sup_{u\in \R}|Q_n(u)|\| = {\cal O}(1)$. By
(\ref{eq:SupHweighted}) of Lemma \ref{lem:supHweighted} and Lemma
\ref{lem:Intsum},
\begin{eqnarray*}
\E\left[\sup_{s\in \R} |Q'_n(s)|^2 \right]
 &\le& C_{\gamma'} \int_\R \|Q'_n(s)\|^2 w_{1+\gamma'} (d s)
 + C_{\gamma'} \int_\R \|Q''_n(s)\|^2 w_{-1 -\gamma'} (d s)\\
 &\le& C_{\gamma'} \sigma^2( f_\varepsilon, w_{1 + \gamma'})
 + C_{\gamma'} \sigma^2( f'_\varepsilon,
   w_{-1 - \gamma'}) = {\cal O}(1).
\end{eqnarray*}
Then there exists $C_1 < \infty$ such that for all $n \in \N$,
\begin{eqnarray*}
\E\left[ \sup_{s, t\in [-r,r], ~ 0\le s-t \le \delta}
|Q_n(s)(1+|s|)^{\gamma / q} - Q_n(t)(1+|t|)^{\gamma / q}|^2
\right] \le \delta^2 C_1.
\end{eqnarray*}
Notice that the upper bound in (\ref{eq:Qnrtail}) goes to $0$ as
$r \to \infty$. Hence (ii) obtains. \qed

\subsection{Proof of Theorem \ref{th:RnW}.}\label{sec:83}
Observe that $(\partial / \partial \theta) {\cal P}_0
F_\varepsilon(\theta |{\cal F}_j) = {\cal P}_0
f_\varepsilon(\theta |{\cal F}_j)$ and ${\cal P}_0
F_\varepsilon(\theta |{\cal F}_j) = 0$ when $\theta = \pm \infty$.
By (\ref{eq:SupHmu1infty}) of Lemma \ref{lem:supHweighted},
\begin{eqnarray*}
\sup_{\theta \in \R} [| {\cal P}_0 F_\varepsilon(\theta |{\cal
F}_j) |^2 ( 1 + |\theta|)^{2 \gamma / q} ]
 \le {q \over {2 \gamma}} \int_\R
 | {\cal P}_0 f_\varepsilon(\theta |{\cal F}_j) |^2
 w_{1 + 2\gamma /q} d \theta.
\end{eqnarray*}
Hence by (\ref{eq:fwep6}),
\begin{eqnarray*}
\sum_{i=0}^\infty \sup_{\theta \in \R} \| {\cal P}_0
 F_\varepsilon(\theta |{\cal F}_j) \|
 \le \sqrt{ q/(2\gamma)}
 \sigma(f_\varepsilon, w_{1 + 2\gamma / q})
 < \infty,
\end{eqnarray*}
which by Lemma \ref{lem:hannan} entails the finite-dimensional
convergence. Since $R_n(s) = G_n(s) + Q_n(s)$, the tightness and
(\ref{eq:Enr}) follows from Lemmas \ref{lem:Gnr} and
\ref{lem:Qntight}. \qed

\subsection{Proof of Theorem \ref{th:modulus}.}\label{sec:84}
Note that $(\log n)^q n^{\max(1, q/4) - q/2} = (\log n)^q n^{1 -
q/2} = {\cal O} (\delta_n^{q/2 - 1})$. By (\ref{eq:ossi-emp-d}) of
Lemma \ref{lem:ossi-emp}, under the proposed condition we have
uniformly in $a$ that
\begin{eqnarray*}
 \E \left[\sup_{0\le s< \delta_n} |G_n(a + s) - G_n(a)|^q \right]
 &\le& C_q (\log^q n) n^{\max(1, q/4) - q/2}
 [ F(a +\delta_n) - F(a)]\cr
 &+& C_q \delta_n^{q/2 - 1} \tau^{ q/2 -1}
 \int_a^{a + \delta_n} f(u) d u \cr
 &\le& C \delta_n^{q/2 - 1} [ F(a +\delta_n) - F(a)].
\end{eqnarray*}
Here the constant $C$ only depends on $\tau, \gamma, q$ and
$\E(|X_1|^\gamma)$. Hence
\begin{eqnarray*}
 &&\sum_{k \in \Z} ( 1 + |k \delta_n|)^{\gamma}
 \E \left[\sup_{0\le s< \delta_n}
 |G_n(k\delta_n + s) - G_n(k \delta_n)|^q \right]\\
 &&\le \sum_{k \in \Z} ( 1 + |k \delta_n|)^{\gamma} C
 \delta_n^{q/2 - 1} [ F(k\delta_n +\delta_n) - F(k\delta_n)]
 \le C \delta_n^{q/2 - 1} \E [( 1 + |X_1|)^\gamma].
\end{eqnarray*}
Let $I_k(\delta) = [k\delta, (1+k) \delta]$. Since $0 < \delta <
1$, we have
\begin{eqnarray*}
{1\over 2} \le {1 \over {1 + \delta}} \le \inf_{t \in I_k(\delta)}
{ {1 + |t|} \over {1 + |k \delta|}} \le \sup_{t \in I_k(\delta)} {
{1 + |t|} \over {1 + |k \delta|}} \le 1 + \delta \le 2
\end{eqnarray*}
and
\begin{eqnarray*}
\sup_{t\in I_k(\delta_n), 0\le s < \delta_n} |G_n(t+s) - G_n(t)|
 &\le& \sup_{t\in I_k(\delta_n), 0\le s < \delta_n}
     |G_n(t+s) - G_n( k \delta_n)| \\
 &&    + \sup_{t\in I_k(\delta_n), 0\le s < \delta_n} |
  G_n( k \delta_n) - G_n(t)| \\
 &\le& \sup_{0\le u < 2\delta_n}
     |G_n(k \delta_n + u) - G_n( k \delta_n)| \\
 &&    + \sup_{0\le s < \delta_n} |
  G_n( k \delta_n + s) - G_n(k \delta_n )| \\
 &\le& 2 \sup_{0\le u < 2\delta_n}
     |G_n(k \delta_n + u) - G_n( k \delta_n)|
\end{eqnarray*}
Therefore,
\begin{eqnarray*}
 &&\E \left[\sup_{t\in \R} ( 1 + |t|)^{\gamma}
 \sup_{0\le s< \delta_n} |G_n(t+s) - G_n(t)|^q \right] \\
 && \le \sum_{k \in \Z} \E
 \left[\sup_{t\in I_k(\delta_n)} ( 1 + |t|)^{\gamma}
 \sup_{0\le s< \delta_n } |G_n(t+s) - G_n(t)|^q \right] \\
 &&\le C \sum_{k \in \Z} ( 1 + |k \delta_n|)^{\gamma}
 \E \left[\sup_{0\le s < 2 \delta_n}
 |G_n(k\delta_n + s) - G_n(k \delta_n)|^q \right]
 \le C \delta_n^{q/2 - 1}.
\end{eqnarray*}
Note that $R_n(s) = G_n(s) + Q_n(s)$. Then (\ref{eq:wmodulus})
follows if it holds with $R_n$ replaced by $G_n$ and $Q_n$
respectively. The former is an easy consequence of the preceding
inequality and Jensen's inequality. To show that
(\ref{eq:wmodulus}) holds with $R_n$ replaced by $Q_n$, let
$\gamma' = 2\gamma / q$. By (\ref{eq:SupHweighted}) of Lemma
\ref{lem:supHweighted} and Lemma \ref{lem:Intsum},
\begin{eqnarray*}
\E \left[ \sup_{x\in \R}( 1 + |x|)^{\gamma'}
 |Q'_n(x)|^2 \right]
 &\le& C \int_\R \|Q'_n(x)\|^2 w_{\gamma' - \mu}(d x)
 + C \int_\R \|Q''_n(x)\|^2 w_{\gamma' + \mu}(d x)\\
 &\le& C\sigma^2 (f_\varepsilon, w_{\gamma'-\mu})
 +C\sigma^2(f'_\varepsilon, w_{\gamma' + \mu})
 < \infty,
\end{eqnarray*}
which entails that
\begin{eqnarray*}
 &&\E \left\{ \sup_{t\in \R} \left[( 1 + |t|)^{\gamma'}
 \sup_{|s| \le \delta_n} |Q_n(t+s) - Q_n(t)|^2 \right] \right\}
 \cr
 &&\le \delta_n^2 \E \left\{\sup_{t\in \R} \left[ ( 1 + |t|)^{\gamma'}
 \sup_{|s| \le \delta_n} |Q'_n(t+s)|^2 \right] \right\}\cr
 &&\le C \delta_n^2 \E \left\{\sup_{x\in \R} [( 1 + |x|)^{\gamma'}
 |Q'_n(x)|^2] \right\} = O(\delta_n^2)
\end{eqnarray*}
and completes the proof. \qed

\begin{remark}{\rm
It is worthwhile to note that the modulus of continuity of $G_n$
has the order $\delta_n^{1 - 2/q}$, while that of $Q_n$ has a
higher order $\delta_n$. \qed}
\end{remark}

\section{Proof of Proposition \ref{prop:Mngtight}}
\label{sec:proofThm2} We shall adopt the truncation technique to
deal with $M_n$. For $r > 0$ define the function $g {\bf 1}_{|
\cdot| > r}$ by $(g {\bf 1}_{| \cdot| > r})(x) = g(x) {\bf 1}_{|x|
> r}$ and
\begin{eqnarray}
\label{eq:Mngr} M_n(g {\bf 1}_{| \cdot| > r}) = n^{-1/2}
\sum_{k=1}^n \{ g(X_k) {\bf 1}_{|X_k| > r} - \E[ g(X_k) {\bf
1}_{|X_k| > r} | {\cal F}_{k-1}] \}.
\end{eqnarray}
The function $g {\bf 1}_{| \cdot| \le r}$ and the process $M_n(g
{\bf 1}_{| \cdot| \le r})$ are similarly defined. Since $M_n(g) =
M_n(g {\bf 1}_{| \cdot| > r}) + M_n(g {\bf 1}_{| \cdot| \le r})$,
the tightness of $\{M_n(g):~ g \in {\cal G}_{\gamma, \mu}\}$
follows from Lemmas \ref{lem:M_outside} and \ref{lem:MnIntight}.
To see this, for any $\delta, \eta > 0$, by Lemma
\ref{lem:M_outside}, there exists $r > 0$ such that
\begin{eqnarray}
\label{eq:PM_outside}  \limsup_{n \to \infty} \p^* \left \{
\sup_{g \in {\cal G}_{\gamma, \mu}} |M_n(g {\bf 1}_{| \cdot| >
r})| \ge {\delta \over 4} \right\} \le {\eta \over 4}.
\end{eqnarray}
By Lemma \ref{lem:MnIntight}, there exists $U_1, \ldots, U_I$ with
$I < \infty$ such that ${\bf 1}_{|\cdot| \le r} {\cal G}_{\gamma,
\mu} \subset \cup_{i=1}^I U_i$ and
\begin{eqnarray}
\label{eq:PM_inside}  \limsup_{n \to \infty} \p^* \left \{
\max_{1\le i\le I} \sup_{g, h \in U_i} |M_n(g-h) | \ge {\delta
\over 4} \right\} \le {\eta \over 4}.
\end{eqnarray}
Let $T_i = \{ h + p {\bf 1}_{|\cdot| > r}:~ h \in U_i, ~ p \in
{\cal G}_{\gamma, \mu}\}$, $1\le i\le I$. Then ${\cal G}_{\gamma,
\mu} \subset \cup_{i=1}^I T_i$ and we have
\begin{eqnarray*}
&&\limsup_{n \to \infty} \p^*
 \left \{ \max_{1\le i\le I} \sup_{g, h \in T_i}
 |M_n(g-h) | \ge \delta \right\} \\
 && \le \limsup_{n \to \infty} \p^*
 \left \{ \max_{1\le i\le I} \sup_{g, h \in T_i}
 |M_n(g{\bf 1}_{|\cdot| \le r} - h{\bf 1}_{|\cdot| \le r}) |
 \ge {\delta \over 2} \right\}\\
 &&+
 \limsup_{n \to \infty} \p^*
 \left \{ \max_{1\le i\le I} \sup_{g, h \in T_i}
 |M_n(g{\bf 1}_{|\cdot| > r}-h{\bf 1}_{|\cdot| > r}) |
 \ge {\delta \over 2} \right\} \\
 && \le \limsup_{n \to \infty} \p^*
 \left \{ \max_{1\le i\le I} \sup_{g, h \in U_i}
 |M_n(g-h) | \ge {\delta \over 2} \right\} \\
 &&+ \limsup_{n \to \infty} \p^*
  \left \{ \sup_{g \in {\cal G}_{\gamma, \mu}}
  |M_n(g {\bf 1}_{| \cdot| > r})|
  \ge {\delta \over 4} \right\}\le {\eta \over 2}
\end{eqnarray*}
in view of (\ref{eq:PM_outside}) and (\ref{eq:PM_inside}). Thus by
definition $\{M_n(g):~ g \in {\cal G}_{\gamma, \mu}\}$ is tight
since $I < \infty$ and $\delta$ and $\eta$ are arbitrarily chosen.

The finite-dimensional convergence is a direct consequence of the
martingale central limit theorem.  The case when ${\cal G} = {\cal
H}_{\eta, \delta}$ can be similarly proved. \qed

\begin{lemma}
\label{lem:M_outside} (i) Assume $\E(|X_1|^\gamma) < \infty$,
$\gamma \ge 0$. Then
\begin{eqnarray}
\label{eq:M_outside} \lim_{t \to \infty} \limsup_{n \to \infty}
\E^* \left \{ \sup_{g \in {\cal G}_{\gamma, \mu}} |M_n(g {\bf
1}_{| \cdot| > t})|^2 \right\} = 0.
\end{eqnarray}
(ii) Under conditions of (ii) of Proposition \ref{prop:Mngtight},
\begin{eqnarray}
\label{eq:M_outsideH} \lim_{t \to \infty} \limsup_{n \to \infty}
\E^* \left \{ \sup_{g \in {\cal H}_{\eta, \delta} } |M_n(g {\bf
1}_{| \cdot| > t})| \right\} = 0.
\end{eqnarray}
\end{lemma}
\noindent{\it Proof.} (i) First assume $\mu < 1$. We shall
generalize the argument in Gin\'e and Zinn (1986). Let $\alpha =
(1- \mu)^{-1}$, $g_r(x) = g(x) {\bf 1}_{x > r^\alpha }$, $r \in
\N$, and the interval $I_j = (j^\alpha, (j + 1 )^\alpha]$. Write
$M_n(g_r) = A_n(g;r) + B_n(g;r)$, where
\begin{eqnarray*}
A_n(g;r) = n^{-1/2} \sum_{j=r}^\infty \sum_{k=1}^n \{ [g(X_k) -
g(j^\alpha)] {\bf 1}_{X_k \in I_j} - \E[ (g(X_k) - g(j^\alpha))
{\bf 1}_{X_k \in I_j} | {\cal F}_{k-1}] \}
\end{eqnarray*}
and
\begin{eqnarray}
\label{eq:Bngr} B_n(g;r) = n^{-1/2} \sum_{j=r}^\infty \sum_{k=1}^n
\{ g(j^\alpha) {\bf 1}_{X_k \in I_j} - \E[ g(j^\alpha) {\bf
1}_{X_k \in I_j} | {\cal F}_{k-1}] \}.
\end{eqnarray}
Let $Z_{j, n} = n^{-1/2} \sum_{k=1}^n [{\bf 1}_{X_k \in I_j} - \E(
{\bf 1}_{X_k \in I_j} | {\cal F}_{k-1}) ]$. Then $\|Z_{j, n}\|^2
\le \p(X_k \in I_j) =: p_j$. By the Cauchy-Schwarz inequality and
(\ref{eq:SupH55}) of Lemma \ref{lem:supHweighted},
\begin{eqnarray*}
\E \left \{\sup_{g \in {\cal G}_{\gamma, \mu}}
 \sum_{j=r}^\infty |g(j^\alpha) Z_{j, n}| \right\}^2
 &\le& \E \left \{\sup_{g \in {\cal G}_{\gamma, \mu}}
 \sum_{j=r}^\infty g^2(j^\alpha) (1+ j^\alpha)^{-\gamma}
 \times \sum_{j=r}^\infty (1+ j^\alpha)^{\gamma} Z^2_{j, n}
 \right \} \cr
 &\le& C_{\gamma, \mu}
 \sum_{j=r}^\infty (1+ j^\alpha)^{\gamma} \E(Z^2_{j, n})
 \le C_{\gamma, \mu} \sum_{j=r}^\infty (1+j^\alpha)^{\gamma} p_j.
\end{eqnarray*}
Since $B_n(g;r) = \sum_{j=r}^\infty g(j^\alpha) Z_{j, n}$ and
$\E[|X_1|^{1/\alpha})^{\alpha \gamma}] = \E(|X_1|^\gamma) <
\infty$,
\begin{eqnarray*}
\limsup_{r \to \infty} \limsup_{n \to \infty}
 \E^*\left[ \sup_{g\in {\cal G}_{\gamma, \mu}} |B_n(g;r)|^2 \right]
 \le C_{\gamma, \mu} \limsup_{r \to \infty}
 \sum_{j=r}^\infty (1+j^\alpha)^{\gamma} p_j = 0.
\end{eqnarray*}
We now deal with $A_n(g;r)$. Clearly (\ref{eq:M_outside}) follows
if the preceding inequality also holds for $A_n(g;r)$. To this
end, let $S_n(u) = n^{-1/2} \sum_{k=1}^n [{\bf 1}_{X_k \in J(u)} -
\E({\bf 1}_{X_k \in J(u)} | {\cal F}_{k-1})]$, where $J(u) = (u, ~
(\lfloor u^{1/\alpha} \rfloor+1)^\alpha]$. Since $ [g(X_k) -
g(j^\alpha)] {\bf 1}_{X_k \in J(u)} = \int_{I_j} g'(u) {\bf
1}_{X_k \in J(u)} d u$, we have
\begin{eqnarray}
\label{eq:Angr} |A_n(g;r)|
 &\le&  n^{-1/2}
 \sum_{j=r}^\infty \left|\sum_{k=1}^n \int_{I_j} g'(u)
 [{\bf 1}_{X_k \in J(u)}
 - \E({\bf 1}_{X_k \in J(u)} | {\cal F}_{k-1})]
 d u \right| \cr
 &\le& \int_{r^\alpha}^\infty |g'(u)| |S_n(u)| d u.
\end{eqnarray}
By the Cauchy-Schwarz inequality, $|A_n(g;r)|^2 \le
\int_{r^\alpha}^\infty |S_n(u)|^2 w_{\gamma - \mu}(d u)$ since $g
\in {\cal G}_{\gamma, \mu}$. So
\begin{eqnarray*}
\lefteqn{ \limsup_{r \to \infty} \limsup_{n \to \infty} \E^*
\left\{ \sup_{g \in {\cal G}_{\gamma, \mu}} |A_n(g;r)|^2 \right \}
 \le \limsup_{r \to \infty} \int_{r^\alpha}^\infty
 \|S_n(u)\|^2 w_{\gamma - \mu}(d u)} \\
 &\le& \limsup_{r\to \infty} \alpha \int_r^\infty
  \p[ t^\alpha < X_1 \le (t+1)^\alpha] (1+t^\alpha)^{\gamma-\mu}
  t^{\alpha -1} d t = 0
\end{eqnarray*}
in view of $\E[ (|X_k|^{1/\alpha})^{ \alpha(\gamma - \mu) + \alpha
-1}] = \E(|X_k|^\gamma)  < \infty$. In the case that $\mu = 1$,
let $I_j = (2^j, 2^{j+1}]$ and $J(u) = (u, 2 u]$. It is easily
seen that the above argument still works.

(ii) Let ${\cal G} = {\cal H}_{\eta, \delta}$. By (\ref{eq:Bngr})
and (\ref{eq:Angr}), if $0 \le \eta - \delta < 1$, then we have
\begin{eqnarray*}
\E \left[ \sup_{g \in {\cal H}_{\eta, \delta} } |B_n(g;r)| \right]
 \le \sum_{j = r}^\infty ( 1 + j^\alpha)^\eta \sqrt {p_j}
\end{eqnarray*}
and
\begin{eqnarray*}
\E \left[ \sup_{g \in {\cal H}_{\eta, \delta} } |A_n(g;r)| \right]
 &\le& \int_{r^\alpha}^\infty (1 + u)^\delta
 \sqrt{\p[|X_1| \in J(u)]} d u\\
 &\le& \int_r^\infty (1 + v^\alpha)^\delta \sqrt{ \p[v^\alpha \le
 |X_1| \le (v+1)^\alpha]} (\alpha v^{\alpha -1}) d v \\
 &\le& C_{\alpha, \delta} \int_r^\infty  v^{\alpha \eta}
 \sqrt{ \p[v^\alpha \le |X_1| \le (v+1)^\alpha]} d v
\end{eqnarray*}
which in view of (\ref{eq:Ginezinn1}) approaches zero if $r \to
\infty$. The case that $\eta- \delta =1$ can be similarly dealt
with. \qed

\bigskip

\begin{lemma}
\label{lem:MnIntight} Let ${\cal G} = {\cal G}_{\gamma, \mu}$ or
${\cal H}_{\eta, \delta}$. Then for any $r > 0$, the process
$\{M_n(g {\bf 1}_{|\cdot| \le r} ):~ g \in {\cal G}\}$ is tight.
\end{lemma}
\noindent{\it Proof.} Consider first ${\cal G} = {\cal G}_{\gamma,
\mu}$. Recall (\ref{eq:Dedd2}) for the definition of the essential
supremum norm $d_2$. Then $d_2(g) \le \| g \|_\infty := \sup_{x\in
\R} |g(x)|$ and
\begin{eqnarray*}
N(u, {\bf 1}_{|\cdot| \le r} {\cal G}_{\gamma, \mu}, d_2) \le N(u,
{\bf 1}_{|\cdot| \le r} {\cal G}_{\gamma, \mu}, \|\cdot
\|_\infty).
\end{eqnarray*}
Let $C_{\gamma, \mu}$ be the constant in (\ref{eq:SupHweighted})
and define the Sobolev class
\begin{eqnarray*}
{\cal S} = \left \{ h: [-r, r] \mapsto \R: \sup_{x \in [-r, r] }
|h(x)|^2 \le C_{\gamma, \mu} (1 + r)^\gamma \mbox{and} \int_{-r}^r
\! |h'(x)|^2 d x \le (1 + r)^{\gamma - \mu} \right\}.
\end{eqnarray*}
Then there exists a constant $C = C(r, \gamma, \mu)$ such that for
every $\epsilon > 0$,
\begin{eqnarray}
\label{eq:Sobolev} \log N( \epsilon, {\cal S}, \|\cdot\|_\infty)
\le {C \over \epsilon}.
\end{eqnarray}
[cf. Birman and Solomjak (1967) or Theorem 2.7.1 in van der Vaart
and Wellner (1996)]. For every $g \in {\bf 1}_{|\cdot| \le r}
{\cal G}_{\gamma, \mu}$, it is easily seen by
(\ref{eq:SupHweighted}) that $\sup_{x \in [-r, r] } |g(x)|^2 \le
C_{\gamma, \mu} (1 + r)^\gamma$ and $\int_{-r}^r |g'(x)|^2 d x \le
(1 + r)^{\gamma - \mu}$. Hence
\begin{eqnarray}\label{eq:GS}
N(u, {\bf 1}_{|\cdot| \le r} {\cal G}_{\gamma, \mu},
\|\cdot\|_\infty) \le N(u, {\cal S}, \|\cdot\|_\infty)
\end{eqnarray}
and consequently
\begin{eqnarray*}
\int_0^1 \sqrt {\log N(u, {\bf 1}_{|\cdot| \le r} {\cal
G}_{\gamma, \mu}, d_2)} d u \le \int_0^1 \sqrt {\log N(u, {\cal
S}, \|\cdot\|_\infty) } d u < \infty.
\end{eqnarray*}
Therefore the lemma follows from Theorem 3.3 in Dedecker and
Louhichi (2002) [see also Section 4.2 therein]. The case that
${\cal G} = {\cal H}_{\eta, \delta}$ can be similarly proved. \qed

\bigskip

{\bf Acknowledgments.} The author is grateful to the referee and
the editor for many helpful comments. The author also thanks
Professors S\'andor Cs\"org\H{o} and Jan Mielniczuk for useful
suggestions.

\bigskip

\centerline{REFERENCES}

\bigskip

\par\noindent\hangindent2.3em\hangafter 1
{\sc Andrews, D. W. K.} and {\sc Pollard, D.} (1994). An
introduction to functional central limit theorems for dependent
stochastic processes. {\it Internat. Statist. Rev.} {\bf 62}
119--132.

\par\noindent\hangindent2.3em\hangafter 1
{\sc Bae, J.} and {\sc Levental, S.} (1995). Uniform CLT for
Markov chains and its invariance principle: A martingale approach.
{\it J. Theoret. Probab.} {\bf 8} 549-570.

\par\noindent\hangindent2.3em\hangafter 1
{\sc Berkes, I.} and {\sc Horv\'ath, L.} (2004). The efficiency of
the estimators of the parameters in GARCH processes. {\it Ann.
Statist.} {\bf 32} 633--655.

\par\noindent\hangindent2.3em\hangafter 1
{\sc Billingsley, P.} (1968). {\it Convergence of Probability
Measures.} New York, Wiley.

\par\noindent\hangindent2.3em\hangafter 1
{\sc Birman, M. S.} and {\sc Solomjak, M. Z.} (1967).
Piecewise-polynomial approximations of functions of the classes
$W_p$. {\it Math. USSR-Sb.} {\bf 73} 295--317.

\par\noindent\hangindent2.3em\hangafter 1
{\sc Bradley, R.}, (2002). {\it Introduction to strong mixing
conditions, Volumes 1 and 2.} Technical report, Indiana
University, Bloomington.

\par\noindent\hangindent2.3em\hangafter 1
{\sc Chow, Y. S.} and {\sc Teicher, H.} (1988). {\it Probability
Theory}, 2nd ed. Springer, New York.

\par\noindent\hangindent2.3em\hangafter 1
{\sc Cs\"org\H o, M.,  Cs\"org\H o, S., Horv\'ath, L.} and {\sc
Mason, D. M.} (1986). Weighted empirical and quantile processes.
{\it Ann. Probab.} {\bf 14} 31-85.

\par\noindent\hangindent2.3em\hangafter 1
{\sc Cs\"org\H{o}, S.} and {\sc Mielniczuk, J.} (1996). The
empirical process of a short-range dependent stationary sequence
under Gaussian subordination. {\it Probab. Theory Related Fields}
{\bf 104} 15--25.

\par\noindent\hangindent2.3em\hangafter 1
{\sc Dedecker, J.} and {\sc Louhichi, S.} (2002). Maximal
inequalities and empirical central limit theorems. In {\it
Empirical Process Techniques for Dependent Data} (Dehling, Mikosch
and Sorensen, editors) 137-159. Birkhauser.

\par\noindent\hangindent2.3em\hangafter 1
{\sc Dedecker, J.} and {\sc Prieur, C.} (2003a). Coupling for
$\tau$-dependent sequences and applications. {\it Tech. Report} \#
2003-2, LSTA, Universit\'e Pierre et Marie Curie-Paris 6.

\par\noindent\hangindent2.3em\hangafter 1
{\sc Dedecker, J.} and {\sc Prieur, C.} (2003b). New dependent
coefficients. Examples and applications to statistics. {\it Tech.
Report} \# 2003-6, LSTA, Universit\'e Pierre et Marie Curie-Paris
6.

\par\noindent\hangindent2.3em\hangafter 1
{\sc Dedecker, J.} and {\sc Rio, E.} (2000). On the functional
central limit theorem for stationary processes. {\it Ann. Inst. H.
Poincar\'e Probab. Statist.} {\bf 36} 1--34.

\par\noindent\hangindent2.3em\hangafter 1
{\sc Dehling, H.}, {\sc Mikosch, T.} and {\sc Sorensen, M.} (eds),
(2002). {\it Empirical Process Techniques for Dependent Data},
Boston: Birkh\"auser.

\par\noindent\hangindent2.3em\hangafter 1
{\sc Dehling, H.} and {\sc Taqqu, M. S.} (1989). The empirical
process of some long-range dependent sequences with an application
to $U$-statistics. {\it Ann. Statist.} {\bf 17} 1767--1786.

\par\noindent\hangindent2.3em\hangafter 1
{\sc Diaconis, P.} and {\sc Freedman, D.} (1999). Iterated random
functions. {\it SIAM Rev.} {\bf 41} 41--76.

\par\noindent\hangindent2.3em\hangafter 1
{\sc Donsker, M. D.} (1952). Justification and extension of Doob's
heuristic approach to the Kolmogorov-Smirnov theorems. {\it Ann.
Math. Statist.} {\bf 23} 277--281.

\par\noindent\hangindent2.3em\hangafter 1
{\sc Doob, J.} (1953). {\it Stochastic Processes}. Wiley.

\par\noindent\hangindent2.3em\hangafter 1
{\sc Douc, R., Fort, G., Moulines, E.} and {\sc Soulier, P.}
(2004). Practical drift conditions for subgeometric rates of
convergence. {\it Ann. Appl. Probab.} {\bf 14} 1353–-1377.

\par\noindent\hangindent2.3em\hangafter 1
{\sc Doukhan, P.} (1994). {\it Mixing. Properties and examples.}
Springer, New York.

\par\noindent\hangindent2.3em\hangafter 1
{\sc Doukhan, P.} (2003). Models, inequalities, and limit theorems
for stationary sequences. In {\it Theory and applications of
long-range dependence} (P. Doukhan, G. Oppenheim and M. S. Taqqu,
eds.) 43--100. Birkh\"auser, Boston, MA.

\par\noindent\hangindent2.3em\hangafter 1
{\sc Doukhan, P., Lang, G.} and {\sc Surgailis, D.} (2002).
Asymptotics of weighted empirical processes of linear fields with
long-range dependence. {\it Ann. Inst. H. Poincar\'e Probab.
Statist.} {\bf 38} 879--896.

\par\noindent\hangindent2.3em\hangafter 1
{\sc Doukhan, P., Lang, G., Surgailis, D.} and {\sc Viano, M. C.}
(2004). Functional limit theorem for the empirical process shifts
with long memory. {\it To appear, J. Theor. Probab.}

\par\noindent\hangindent2.3em\hangafter 1
{\sc Doukhan, P.} and {\sc Louhichi, S.} (1999). A new weak
dependence condition and applications to moment inequalities. {\it
Stochastic Process. Appl.} {\bf 84} 313-342.

\par\noindent\hangindent2.3em\hangafter 1
{\sc Doukhan, P., Massart, P.} and {\sc Rio, E.} (1995).
Invariance principles for absolutely regular empirical processes.
{\it Ann. Inst. H. Poincar\'e Probab. Statist.} {\bf 31} 393--427.

\par\noindent\hangindent2.3em\hangafter 1
{\sc Doukhan, P.} and {\sc Surgailis, D.} (1998). Functional
central limit theorem for the empirical process of short memory
linear processes. {\it C. R. Acad. Sci. Paris Ser. I Math. } {\bf
326} 87--92.

\par\noindent\hangindent2.3em\hangafter 1
{\sc Dudley, R. M.} (1978). Central limit theorems for empirical
measures. {\it Ann. Probab.} {\bf 6}, 899--929.

\par\noindent\hangindent2.3em\hangafter 1
{\sc Einmahl, J. H. J.} and {\sc Mason, D. M.} (1988). Strong
limit theorems for weighted quantile processes. {\it Ann. Probab.}
{\bf 16} 1623-1643.

\par\noindent\hangindent2.3em\hangafter 1
{\sc Elton, J. H.} (1990). A multiplicative ergodic theorem for
Lipschitz maps. {\it Stochastic Process. Appl.} {\bf 34} 39--47.

\par\noindent\hangindent2.3em\hangafter 1
{\sc Feller, W.} (1971). {\it An introduction to probability
theory and its applications.} Vol. II. John Wiley \& Sons, New
York.

\par\noindent\hangindent2.3em\hangafter 1
{\sc Gastwirth, J. L.} and {\sc Rubin, H.} (1975). The asymptotic
distribution theory of the empiric ${\rm cdf}$ for mixing
stochastic processes. {\it Ann. Statist.} {\bf 3} 809--824.

\par\noindent\hangindent2.3em\hangafter 1
{\sc Gin\'e, E.} and {\sc Zinn, J} (1984). Some limit theorems for
empirical processes (in Special Invited Papers). {\it Ann.
Probab.} {\bf 12} 929-989.

\par\noindent\hangindent2.3em\hangafter 1
{\sc Gin\'e, E.} and {\sc Zinn, J} (1986). Empirical processes
indexed by Lipschitz functions. {\it Ann. Probab.} {\bf 14}
1329-1338.

\par\noindent\hangindent2.3em\hangafter 1
{\sc Giraitis, L.} and {\sc Surgailis, D.} (1999). Central limit
theorem for the empirical process of a linear sequence with long
memory. {\it J. Statist. Plann. Inference} {\bf 80} 81--93.

\par\noindent\hangindent2.3em\hangafter 1
{\sc Gordin, M. I.} (1969). The central limit theorem for
stationary processes. {\it Dokl. Akad. Nauk SSSR} {\bf 188}
739--741.

\par\noindent\hangindent2.3em\hangafter 1
{\sc Gordin, M. I.} and {\sc Lifsic, B.} (1978). The central limit
theorem for stationary Markov processes. {\it  Dokl. Akad. Nauk
SSSR} {\bf 239} 766--767.

\par\noindent\hangindent2.3em\hangafter 1
{\sc Gordin, M. I.} and {\sc Holzmann, H.} (2004). The central
limit theorem for stationary Markov chains under invariant
splittings. {\it Stoch. Dynamics} {\bf 4} 15--30.

\par\noindent\hangindent2.3em\hangafter 1
{\sc Hall, P.} and {\sc Heyde, C. C.} (1980). {\it Martingale
limit theory and its applications}. Academic Press, New York.

\par\noindent\hangindent2.3em\hangafter 1
{\sc Hannan, E. J.} (1973). Central limit theorems for time series
regression. {\it Z. Wahrsch. und Verw. Gebiete} {\bf 26} 157--170.

\par\noindent\hangindent2.3em\hangafter 1
{\sc Ho, H. C.} and {\sc Hsing, T.} (1996). On the asymptotic
expansion of the empirical process of long-memory moving averages.
{\it Ann. Statist.} {\bf 24} 992--1024.

\par\noindent\hangindent1.5em\hangafter 1
{\sc Ho, H. C.} and {\sc Hsing, T.} (1997). Limit theorems for
functionals of moving averages {\it Ann. Probab.} {\bf 25}
1636-1669.

\par\noindent\hangindent2.3em\hangafter 1
{\sc Hsing, T.} (1999). On the asymptotic distributions of partial
sums of functionals of in.nitevariance moving averages. {\it Ann.
Probab.} {\bf 27} 1579--1599.

\par\noindent\hangindent2.3em\hangafter 1
{\sc Hsing, T.} and {\sc Wu, W. B.} (2004). On weighted
$U$-statistics for stationary processes. {\it Ann. Probab.} {\bf
32} 1600--1631.

\par\noindent\hangindent2.3em\hangafter 1
{\sc Ibragimov, I. A.} (1962). Some limit theorems for stationary
processes. {\it Theory Probab. Appl.} {\bf 7} 349-382.

\par\noindent\hangindent2.3em\hangafter 1
{\sc Jarner, S.} and {\sc Tweedie, R.} (2001). Locally contracting
iterated random functions and stability of Markov chains. {\it J.
Appl. Probab.} {\bf 38} 494--507.

\par\noindent\hangindent2.3em\hangafter 1
{\sc Kipnis, C.} and {\sc Varadhan, S. R. S.} (1986). Central
limit theorem for additive functionals of reversible Markov
processes and applications to simple exclusions. {\it Comm. Math.
Phys.} {\bf 104} 1--19.

\par\noindent\hangindent2.3em\hangafter 1
{\sc Koul, H.} and {\sc Surgailis, D.} (2001). Asymptotics of
empirical processes of long memory moving averages with infinite
variance. {\it Stochastic Process. Appl.} {\bf 91} 309-336.

\par\noindent\hangindent2.3em\hangafter 1
{\sc Major, P.} (1981). {\it Multiple Wiener-Ito integrals: with
applications to limit theorems}. Berlin; New York: Springer.

\par\noindent\hangindent2.3em\hangafter 1
{\sc Mehra, K. L.} and {\sc Rao, M. S.} (1975) Weak convergence of
generalized empirical processes relative to $d\sb{q}$ under strong
mixing. {\it Ann. Probab.} {\bf 3} 979--991.

\par\noindent\hangindent2.3em\hangafter 1
{\sc Nishiyama, Y.} (2000). Weak convergence of some classes of
martingales with jumps. {\it Ann. Probab.} {\bf 28} 685--712.

\par\noindent\hangindent2.3em\hangafter 1
{\sc Opic, B.} and {Kufner, A.} (1990). {\it Hardy-type
inequalities.} Longman Scientific \& Technical; New York; Wiley.

\par\noindent\hangindent2.3em\hangafter 1
{\sc Ossiander, M.} (1987). A central limit theorem under metric
entropy with $L\sb 2$ bracketing. {\it Ann. Probab.} {\bf 15}
897--919.

\par\noindent\hangindent2.3em\hangafter 1
{\sc Peligrad, M.} (1996). On the asymptotic normality of
sequences of weak dependent random variables. {\it J. Theor.
Probab.}, {\bf 9} 703-715.

\par\noindent\hangindent2.3em\hangafter 1
{\sc Pipiras, V.} and {\sc Taqqu, M. S.} (2003). Central limit
theorems for partial sums of bounded functionals of
infinite-variance moving averages. {\it Bernoulli} {\bf 9}
833--855.

\par\noindent\hangindent2.3em\hangafter 1
{\sc Pollard, D.} (1984). {\it Convergence of stochastic
processes}. Springer, New York.

\par\noindent\hangindent2.3em\hangafter 1
{\sc Pollard, D.} (2002). Maximal inequalities via bracketing with
adaptive truncation. {\it Ann. Inst. H. Poincar\'e Probab.
Statist.} {\bf 38} 1039--1052.

\par\noindent\hangindent2.3em\hangafter 1
{\sc Prieur, C.} (2002). An empirical functional central limit
theorem for weakly dependent sequences. {\it Probab. Math.
Statist.} {\bf 22} 259--287.

\par\noindent\hangindent2.3em\hangafter 1
{\sc Rio, E.} (1998). Processus empiriques absolument r\'eguliers
et entropie universelle. {\it Probab. Theory Related Fields} {\bf
111} 585--608.

\par\noindent\hangindent2.3em\hangafter 1
{\sc Rio, E.} (2000). {\it Theorie asymptotique des processus
aleatoires faiblement dependants}. Math\'ematiques et Applications
{\bf 31}. Springer, Berlin

\par\noindent\hangindent2.3em\hangafter 1
{\sc Rosenblatt, M.} (1956). {A central limit theorem and a strong
mixing condition.} {\it Proc. Nat. Acad. Sci. USA} {\bf 42}
43--47.

\par\noindent\hangindent2.3em\hangafter 1
{\sc Shao, Q. M.} and {\sc Yu, H.} (1996). Weak convergence for
weighted empirical processes of dependent sequences. {\it Ann.
Probab.} {\bf 24} 2098-2127.

\par\noindent\hangindent2.3em\hangafter 1
{\sc Shorack, G. R.} and {\sc Wellner, J. A.} (1986). {\it
Empirical processes with applications to statistics.} John Wiley
\& Sons, New York.

\par\noindent\hangindent2.3em\hangafter 1
{\sc Steinsaltz, D.} (1999). Locally contractive iterated function
systems. {\it Ann. Probab.} {\bf 27} 1952--1979.

\par\noindent\hangindent2.3em\hangafter 1
{\sc Straumann, D.} and {\sc Mikosch, T.} (2003). Quasi-MLE in
heteroscedastic times series: a stochastic recurrence equations
approach. {\it Technical Report},  Institute for Mathematical
Sciences, University of Copenhagen.

\par\noindent\hangindent2.3em\hangafter 1
{\sc Surgailis, D.} (2002). Stable limits of empirical processes
of moving averages with infinite variance. {\it Stochastic
Process. Appl.} {\bf 100} 255-274.

\par\noindent\hangindent2.3em\hangafter 1
{\sc Van der Vaart, A. W.} and {\sc Wellner, J. A.} (1996). {\it
Weak convergence and empirical processes}. Springer-Verlag, New
York.

\par\noindent\hangindent2.3em\hangafter 1
{\sc Van der Vaart, A. W.} (1996). New Donsker classes. {\it Ann.
Probab.} {\bf 24} 2128--2140.

\par\noindent\hangindent2.3em\hangafter 1
{\sc Voln\'y, D.} (1993). Approximating martingales and the
central limit theorem for strictly stationary processes. {\it
Stochastic Process. Appl.} {\bf 44} 41--74.

\par\noindent\hangindent2.3em\hangafter 1
{\sc Withers, C. S.} (1975). Convergence of empirical processes of
mixing rv's on $[0,\,1]$. {\it Ann. Statist.} {\bf 3} 1101--1108.

\par\noindent\hangindent2.3em\hangafter 1
{\sc Woodroofe, M.} (1992). A central limit theorem for functions
of a Markov chain with applications to shifts. {\it Stochastic
Process. Appl.} {\bf 41} 33--44.

\par\noindent\hangindent2.3em\hangafter 1
{\sc Wu, W. B.} (2003a). Empirical processes of long-memory
sequences. {\it Bernoulli} {\bf 9} 809--831.

\par\noindent\hangindent2.3em\hangafter 1
{\sc Wu, W. B.} (2003b). Additive functionals of infinite-variance
moving averages. {\it Statistica Sinica} {\bf 13} 1259--1267

\par\noindent\hangindent2.3em\hangafter 1
{\sc Wu, W. B.} (2004a). On the Bahadur representation of sample
quantiles for stationary sequences. {\it To appear, Ann. Statist.}

\par\noindent\hangindent2.3em\hangafter 1
{\sc Wu, W. B.} (2004b). Fourier transforms of stationary
processes {\it To appear, Proc. Amer. Math. Soc.}

\par\noindent\hangindent2.3em\hangafter 1
{\sc Wu, W. B.} and {\sc Mielniczuk, J.} (2002). Kernel density
estimation for linear processes. {\it Ann. Statist.} {\bf 30}
1441--1459.

\par\noindent\hangindent2.3em\hangafter 1
{\sc Wu, W. B.} and {\sc Woodroofe, M.}, (2000). A central limit
theorem for iterated random functions. {\it J. Appl. Probab.} {\bf
37} 748--755.

\par\noindent\hangindent2.3em\hangafter 1
{\sc Wu, W. B.} and {\sc Woodroofe, M.}, (2004). Martingale
approximations for sums of stationary processes. {\it Ann.
Probab.} {\bf 32} 1674--1690.

\par\noindent\hangindent2.3em\hangafter 1
{\sc Wu, W. B.} and {\sc Shao, X.} (2004). Limit theorems for
iterated random functions. {\it J. Appl. Probab.} {\bf 41}
425--436.

\smallskip

\noindent{\n Department of Statistics,}\\
\noindent{\n The University of Chicago,}\\
\noindent{\n 5734 S. University Avenue, Chicago, IL 60637}\\
{\tt wbwu@galton.uchicago.edu}

\end{document}